# Multi-valued Logics

Giovanni Panti
Dept. of Mathematics and Computer Science
via delle Scienze, 208
33100 Udine, Italy
panti@dimi.uniud.it

# Contents

# 1 Introduction

It is almost a safe bet to assume that most students, after receiving the first rudiments of Boolean calculus, wonder: "Why not more than two truth-values?" This question is both natural and vague; attempts to give an answer originated a lot of logical systems. Because of the original vagueness, it is not surprising that the answers range over a wide spectrum of possibilities, and that the borders of the discipline are not clear. Up to now, no many-valued system has imposed as *the* system of many-valued logic. The situation is analogous to that occurring in —say— modal logic, and the sizes of the relative literatures can be compared.

## 1.1 What is many-valued logic?

The simplest and probably the best answer is the historical one: many-valued logic is the complex of studies that originated from the papers of Łukasiewicz and Post in the twenties. The idea underlying these studies is to extend the scope of classical logic by considering a set of truth-values larger than the usual $\{0,1\}$. The new set may be finite or infinite and, in most cases, it will bear some order structure, making it a poset, or a lattice, or a chain, with a top element ("complete truth"), and a bottom one ("complete falsity").

Some more conditions are generally assumed:

1. there is a finite number of connectives, each one of finite arity;
2. the connectives are truth-functional;
3. connectives and truth-values have a "logical meaning";
4. fuzziness phenomena are not present at the metalogical level.

While the first condition does not require any particular comment, the second one restricts greatly the possible interpretations of many-valued systems. For example, probabilistic interpretations are ruled out: we may know that the probabilities of the events $A$ and $B$ are both $1/2$, without being able to compute, in absence of information about the stochastic dependence of $A$ and $B$, the probability of ($A$ and $B$). In a similar fashion, an interpretation of the truth-values in terms of modalities is questionable. Take as an example Łukasiewicz's original interpretation of his three-valued logic (see Example 2.1.2): we have the truth-values $0, 1/2, 1$, with $1/2$ to be read as "possible", and a binary connective $\rightarrow$, to be read as "implies". Let us say that $A$ is the proposition "Tomorrow it will rain", while $B$ is the proposition "Riemann's hypothesis will be proved in the next ten years". In Łukasiewicz's interpretation, it is reasonable to assign the value $1/2$ to both $A$ and $B$. Now, most of us will agree to assign the value $1$ to both $A \rightarrow A$ and $B \rightarrow B$, but surely we will feel uncomfortable in assigning $1$ to $A \rightarrow B$, as truth-functionality requires.

The third condition in our list is deliberately vague. Actually, even during the first pioneering period of many-valued logic, when philosophical considerations lead the way, Post systems stood as examples of systems without such a "logical meaning". Today, the debate about the "meaning" is more relaxed, and the stress is about the usefulness of this or that system. Also, there is a tendency not to consider a single system, but rather to formulate general theories about classes of systems. In our opinion, the relevance to be given to the third condition depends very much on one's own tastes. We only remark that many-valued systems have been introduced as extensions of classical, two-valued, logic. Most of them have values to represent "complete truth" and "complete falsity" and, if we restrict to these values, we get two-valued logic. If the truth-values are ordered in some way, and the connnectives are intended to parallel those of two-valued logic, then they must behave in a reasonable way. So, a "conjunction" connective will be isotone in both variables, an "implication" one will be antitone in the first and isotone in the second variable, a "negation" will be order-reversing. One may well let condition (3) drop altogether, but the price to be paid is a loss of some of the characteristics peculiar to the discipline, and a shift towards universal algebra and combinatorics. We refer to [81, §2.18] for a discussion about this issue and further references.

Our fourth condition is intended to establish a distinction between many-valued logic and fuzzy logic. As both of these logics have —appropriately— fuzzy boundaries, this is not easily accomplished, but we put the question in the following terms. Fuzzy logic in the broad sense is everything concerning fuzziness. In the narrow sense, it is the formal logical calculus of fuzziness. We have to distinguish between the cases where fuzziness is present at the logical level (the level of the calculus and its semantics) from the cases where fuzziness is present at the metalogical level (the level of the construction of the calculus and its semantics). For example, we may have a fuzzy statement $A$ which is an axiom of the theory $T$ in a crisp way (i.e., $A$ belongs to the axiom set of $T$ with degree of membership 1), or a crisp statement $B$ which is an axiom of $T$ in a fuzzy way (i.e., $B$ is an axiom with degree of membership $\neq 0, 1$), or any one of the possible combinations. In other words, we distinguish between *truth* and *belief*. Think of some person showing us a picture: the picture may be out of focus, and the person's hand may be trembling. Every combination of effects may take place, and it is important to recognize from where the out-of-focus effect originates. In many-valued logic, the hand is always steady: no fuzziness phenomenon is to be expected at the metalogical level, and a set is fuzzy only if it is the realization of a predicate symbol in a structure (see §2.6). We will return to this and related issues in the last section.

## 1.2 How this survey is organized

In recent years, due to applications in computer science, many-valued logic attracted renewed interest. The literature, which was already large, has grown considerably. Many-valued systems are now being applied in fields like circuit verification [42], analysis of communication with feedback [67], [68], logic programming [28], [20], [29], artificial intelligence [33], program verification [9], non-monotonic reasoning [94], [23], [91], natural language processing [8], [88].

We have decided to give the basic definitions at the propositional level, and to sketch the extension to the predicate case in a separate section (§2.6). As a matter of fact, we feel that the main ideas and problems can be better understood at the propositional level. Moreover, we have completely omitted some issues that are of specialized interest and are conveniently treated in the literature. In particular, we have omitted the use of many-valued logic to prove the independence of axioms in propositional two-valued logic [17], [36], and the attempts towards a many-valued formulation of the comprehension axiom scheme in set theory (see [104, §2.14], [113] and references therein). Also, we do not discuss many-valued modal logic; see [30], [31], [102], [100], [73], [86], [63].

We develop the basic semantics in §2.1, and we present some proof theory of finite-valued systems in §2.2. We assign a section to each one of the two most studied systems: the Post system (§2.3), and the Łukasiewicz one (§2.4). §2.5 is devoted to a very short selection of some of the other systems; see [83], [81], [12], [36], [57] for many more examples.

The third part is devoted to the algebraic theory of MV-algebras, and is more mathematically oriented; people not interested in MV-algebras may just skip it. Conversely, one can read the third part just after having read sections §2.1 and §2.4.

In the last section we sum up our discussion, and try to give motivations for many-valued logic.

The References section is by no means exhaustive; we refer to [81], [111], [110], [12], [36], [57] for comprehensive bibliographies.

I would like to thank most sincerely Reiner Hähnle, Petr Hájek, Daniele Mundici, Vilém Novák, Zbigniew Stachniak, and Richard Zach, who read previous versions of this survey and provided extremely valuable comments and very constructive and detailed criticism.

## 2 General theory

A language $L$ for a propositional many-valued logic is given by:

1. a finite or denumerable set $p, q, r, \ldots$ of propositional symbols;

2. a finite set of connectives $c_1, \ldots, c_m$, with $c_i$ having arity $u_i \geq 0$.

The set $FORM(L)$ of formulas in the language $L$ is defined as usual: any propositional symbol is a formula and, if $c$ is a connective of arity $u$ and $A_1, \ldots, A_u$ are formulas, then $cA_1 \ldots A_u$ is a formula.

In specific examples, we always have connectives of arity at most 2, and we write binary connectives using infix notation. We also use left and right parentheses according to the standard conventions; so, e.g., $\to p \to qp$ is written $p \to (q \to p)$. We use $A, B, C, \ldots$ to denote formulas, and $\Gamma, \Delta, \ldots$ to denote sets of formulas.

A *logical matrix* for $L$ is a system $\mathcal{M} = (M, D)$, such that $M$ is an abstract algebra of type $(u_1, \ldots, u_m)$ (i.e., a set $M$ together with a collection of functions $c_1^M, \ldots, c_m^M$, with $c_i^M : M^{u_i} \to M$), and $D$ is a subset of $M$. The elements of $M$ are to be thought of as truth-values, while $D$ is the set of *designated* truth-values. 0ary connectives are regarded as logical constants, i.e., elements of $M$. By a standard abuse of notation, we drop the superscript in referring to $c_1^M, \ldots, c_m^M$; in other words, we use the same symbol to refer both to a connective and to the function on the truth-values associated to that connective. We always assume that $M$ has cardinality at least 2, and we say that the pair $(L, \mathcal{M})$ constitutes a propositional many-valued logic.

### 2.1 Semantics

A *structure* for $(L, \mathcal{M})$ is a mapping $\sigma$ from the set of propositional symbols to the algebra $M$. Any structure extends to a homomorphism from the absolutely free algebra $FORM(L)$ to $M$ in the standard way.

If $\Gamma$ is a subset of $FORM(L)$ and $A$ is a formula, then we say that $A$ is a *semantical consequence of* $\Gamma$ *with respect to* $\mathcal{M}$, and we write $\Gamma \models_\mathcal{M} A$, if the following holds:

for any structure $\sigma$ for $(L, \mathcal{M})$, if $\sigma(B) \in D$ for every $B \in \Gamma$, then $\sigma(A) \in D$.

To simplify notation, we write $B_1, \ldots, B_r \models_\mathcal{M} A$ for $\{B_1, \ldots, B_r\} \models_\mathcal{M} A$. We say that $A$ is *valid* if $\models_\mathcal{M} A$, and we say that $A$ is *satisfiable* if there exists a structure $\sigma$ with $\sigma(A) \in D$.

If the matrix $\mathcal{M}$ is finite (i.e., $M$ is finite, say of cardinality $n$), then we identify $M$ either with the set $\mathbf{n} = \{0, 1, \ldots, n-1\}$, or with the set $I_n = \{0, 1/(n-1), \ldots, (n-2)/(n-1), 1\}$, depending on the context or the tradition.

**Example 2.1.1** *Classical logic.* Consider the matrix $\mathcal{B} = (\mathbf{2}, \{1\})$, where $\mathbf{2} = \{0, 1\}$ is the two-element Boolean algebra in the language $(\vee, \neg)$. Then the relation $\models_\mathcal{B}$ is just ordinary Boolean semantical consequence, and the set of valid formulas is the set of Boolean tautologies.

**Example 2.1.2** *Łukasiewicz's three-valued logic.* We have a binary connective $\to$ (implication), and a unary one $\neg$ (negation). The matrix $\mathcal{L}_3 = (I_3, \{1\})$ is defined as follows:

| | $\neg$ |
|---|---|
| 0 | 1 |
| 1/2 | 1/2 |
| 1 | 0 |

| $\to$ | 0 | 1/2 | 1 |
|---|---|---|---|
| 0 | 1 | 1 | 1 |
| 1/2 | 1/2 | 1 | 1 |
| 1 | 0 | 1/2 | 1 |

We consider $I_3$ to be totally ordered ($0 < 1/2 < 1$), with 0 representing "complete falsity" and 1 "complete truth". Observe that $\to, \neg$, if restricted to $\{0, 1\}$, behave like classical implication and negation. It follows that the algebra $I_3$ is not functionally complete; for example, no term built up from $\to, \neg$ may represent the function $: I_3 \to I_3$ whose value is constantly $1/2$.

**Definition 2.1.3** A *consequence operation* is a function $C$ from the powerset of $FORM(L)$ to itself such that:

1. $\Gamma \subseteq C(\Gamma)$;
2. $C(C(\Gamma)) = C(\Gamma)$;
3. $\Gamma \subseteq \Delta$ implies $C(\Gamma) \subseteq C(\Delta)$.

$\Gamma$ is *consistent* if $C(\Gamma) \neq FORM(L)$. $C$ is:

- *structural* if, for every substitution $\alpha$ of formulas for variables, we have $\alpha(C(\Gamma)) \subseteq C(\alpha(\Gamma))$;

- *compact* if $A \in C(\Gamma)$ implies that there exists a finite subset $\Gamma_0$ of $\Gamma$ such that $A \in C(\Gamma_0)$;

- *uniform* if the following holds:

    if $A \in C(\Gamma \cup \Delta)$, $\Gamma$ is consistent, and $\Gamma$ does not share variables with $\Delta \cup \{A\}$, then $A \in C(\Delta)$.

We defined a logic to be a pair $(L, \mathcal{M})$. A more general approach would be to define a logic to be a pair $(L, C)$, where $L$ is a language, and $C$ is a consequence operation on $FORM(L)$. As a matter of fact, for any matrix $\mathcal{M}$, the mapping $C_\mathcal{M}$ defined by

$$C_\mathcal{M}(\Gamma) = \{A \in FORM(L) : \Gamma \models_\mathcal{M} A\}$$

is a consequence operation. Let $C$ be an arbitrary consequence operation. We say that $\mathcal{M}$ is *weakly adequate* for $C$ if $C(\emptyset) = C_\mathcal{M}(\emptyset)$, and say that $\mathcal{M}$ is *strongly adequate* for $C$ if, for every $\Gamma$, $C(\Gamma) = C_\mathcal{M}(\Gamma)$. The consequence operations for which there exists a strongly adequate matrix are characterized in [54], [108]. Note that these strongly adequate matrices may well be uncountable (see §2.5.2). The book [109] is a basic reference for the theory of consequence operations.

## 2.2 Some logical calculi

Suppose we are given a logic $(L, \mathcal{M})$. If the set of valid formulas is recursively enumerable, it is natural to ask for a good calculus. Depending on one's own tastes and needs, good may mean elegant, efficient, easily implemented, cut-free, with cut elimination, or what else. If the set of valid formulas is recursive, then one looks for a good decision procedure. Of course, if the matrix is finite, then truth-tables yield a decision procedure, which is generally unsatisfactory, due to computational costs. The first many-valued calculi were Hilbert-style. There was always an implication connective around, and the only rule of inference was modus ponens. Due to the obvious computational drawbacks of this approach, analytic systems were developed (tableau and sequent systems), where the matrix semantics is closer and the distinction between designated and undesignated values loses its importance.

It is significant that the development and analysis of finite-valued calculi is leading to a clarification and new insight into the structure of classical two-valued calculi. As an example of this feature, we cite the unifying approach to sequent, natural deduction, and tableau calculi developed by Baaz, Fermüller, and Zach [2], [3], [114].

Practically every calculus for 2-valued logic has been extended to many-valued systems, at least for specific logics. We have various forms of Hilbert-style calculi, sequent calculi, natural deduction, resolution, both at the propositional and at the predicate level. Recent proof systems tend not to concentrate on a single logic, but rather to be flexible enough to be applied to a wide range of logics. We will sketch a few systems of this sort, namely sequent and tableau calculi, and resolution. All of them are general-purpose calculi, i.e., if specialized to any finite matrix, they give a calculus which is sound and complete with respect to that matrix. If the matrix is not finite, then the axiomatization is more complex and Hilbert-style systems regain importance. In §3 we will address this issue in one specific case, namely the infinite-valued logic of Łukasiewicz.

We first need some definitions: assume we are given a fixed finite matrix $\mathcal{M}$ of cardinality $n$. A *signed formula* is a formula in the language $L$, labelled with a truth-value; for example, if $\to$ is a binary connective in $L$ and $M = \{0, 1, 2, 3\}$, then $(p \to q)^0$, $(p \to q)^2$, $q^1$ are distinct signed formulas. A *literal* is a signed formula in which the labelled formula is a propositional symbol. A *signed formula expression* is a Boolean combination, with connectives $\neg, \vee, \wedge$, of signed formulas. Any structure $\sigma$ for $(L, \mathcal{M})$ induces a mapping $\overline{\sigma}$ from the set of signed formulas to the two-element Boolean algebra, by naturally setting $\overline{\sigma}(A^i) = 1$ if and only if $\sigma(A) = i$. The mapping $\overline{\sigma}$ extends to the set of all signed formula expressions in the obvious way. We say that the signed formula expressions $E, F$ are *equivalent*, and we write $E \equiv F$, if, for any $\sigma$, we have $\overline{\sigma}(E) = \overline{\sigma}(F)$. Of course, if $E, F$ are equivalent as Boolean formulas, then $E \equiv F$.

As, for any formula $A$ and any truth-value $i$, $\neg(A^i) \equiv \bigvee_{j \neq i} A^j$, it follows that, for any signed formula expression $E$, we may find an equivalent signed formula expression $C$ which is in conjunctive normal form (CNF) and does not contain the connective $\neg$. By applying the distributive laws to $C$, we obtain an equivalent signed formula expression $D$ which is in disjunctive normal form (DNF) and does not contain $\neg$. Note that, for a given $E$, neither $C$ nor $D$ are unique; we refer to [114, §1.5], [40, §4.5] for a discussion about techniques to minimize their length. Note also that we tacitly apply the commutative and associative laws for $\vee, \wedge$ to signed formula expressions in CNF and DNF; in other words, we assume that we are actually working with sets of clauses.

Let $c$ be a connective in $L$ of arity $u$, let $p_1, \ldots, p_u$ be propositional symbols, and let $i$ be a fixed truth-value. By examining the truth-table for $c$, we may easily construct a signed formula expression $T$, containing only signed formulas from the set $\{p_r^j : 1 \leq r \leq u, j \in M\}$, and such that $(cp_1 \ldots p_u)^i \equiv T$. By the above observations, we may construct two (non-unique) signed formula expressions $C, D$ such that:

1. $C, D$ do not contain $\neg$;

2. $C, D$ contain only literals from the set $\{p_r^j : 1 \leq r \leq u, j \in M\}$;

3. $C$ is in CNF and $D$ is in DNF;

4. $T, C, D$ are equivalent as Boolean formulas, and hence they are all equivalent to $(cp_1 \ldots p_u)^i$.

We call $C$ an $i$th CNF for $c$, and $D$ an $i$th DNF for $c$.

**Example 2.2.1** Let $c$ be implication in Łukasiewicz's three-valued logic. Then:

$(p \to q)^0 \equiv p^1 \wedge q^0;$

$(p \to q)^{1/2} \equiv (p^{1/2} \wedge q^0) \vee (p^1 \wedge q^{1/2}) \equiv (p^{1/2} \vee p^1) \wedge (p^{1/2} \vee q^{1/2}) \wedge (q^0 \vee p^1) \wedge (q^0 \vee q^{1/2}) \equiv (p^{1/2} \vee q^{1/2}) \wedge (q^0 \vee p^1)$;

$(p \to q)^1 \equiv p^0 \vee q^1 \vee (p^{1/2} \wedge q^{1/2}) \equiv (p^0 \vee q^1 \vee p^{1/2}) \wedge (p^0 \vee q^1 \vee q^{1/2}).$

### 2.2.1 Sequent and tableau calculi

**Sequent calculi**

Let us fix a finite matrix $\mathcal{M} = (M, D)$ of cardinality $n$; without loss of generality, we denote the elements of $M$ by $0, 1, \ldots, n-1$. A *sequent* is a finite set of signed formulas; we use $G, F, H$ to denote sequents. We write $G, F$ for $G \cup F$ and $G, A^i$ for $G \cup \{A^i\}$. We sometimes write a sequent $G$ in the form $\Gamma_0 \Rightarrow \Gamma_1 \Rightarrow \cdots \Rightarrow \Gamma_{n-1}$, where
$\Gamma_i = \{A : A^i \in G\}$; we use this notation mainly if we think of the truth-values as a totally ordered set, with 0 standing for the falsest value, and $n-1$ for the truest.

A structure $\sigma$ satisfies $G$ if $\overline{\sigma}(\bigvee G) = 1$; in other words, if there exist $1 \leq i \leq n-1$ and $A \in \Gamma_i$ such that $\sigma(A) = i$. A sequent is *valid* if it is satisfied by every structure. If $\Gamma$ is a finite set of formulas, and $A$ is a formula, then $\Gamma \models_\mathcal{M} A$ iff the sequent $\Sigma_0 \Rightarrow \Sigma_1 \Rightarrow \cdots \Rightarrow \Sigma_{n-1}$, where $\Sigma_i = A$ if $i \in D$ and $\Sigma_i = \Gamma$ if $i \notin D$, is valid. In particular, $A$ is valid iff the sequent $\{A^i : i \in D\}$ is valid. It is clear that the definitions of satisfiability and validity for an ordinary two-valued sequent fall into this framework.

Let $c$ be a connective of arity $u$, and let $i$ be a truth-value. Choose your favorite $i$th CNF for $c$, say $(cp_1 \ldots p_u)^i \equiv D_1 \wedge D_2 \wedge \cdots \wedge D_h$, where each $D_s$ is a disjunction of literals from the set $\{p_r^j : 1 \leq r \leq u, j \in M\}$. Then $D_1 \wedge D_2 \wedge \cdots \wedge D_h$ is naturally associated with the following *introduction rule for $c$ at place $i$*:

$$\frac{G, F_1 \quad G, F_2 \quad \cdots \quad G, F_h}{G, (cA_1 \ldots A_u)^i}$$

where $G$ is an arbitrary sequent, and $F_s$ is the sequent obtained from the disjunction $D_s$ (thought of as the set of its disjuncts) by substituting $A_r^j$ for $p_r^j$, for all $1 \leq r \leq u$ and $j \in M$.

**Example 2.2.2** The 1/2th CNF for $\to$ in Łukasiewicz's three-valued logic in Example 2.2.1 is associated with the rule

$$\frac{\Gamma \Rightarrow \Delta, A, B \Rightarrow \Sigma \quad \Gamma, B \Rightarrow \Delta \Rightarrow \Sigma, A}{\Gamma \Rightarrow \Delta, A \to B \Rightarrow \Sigma}$$

**Definition 2.2.3** A *sequent calculus* for a propositional many-valued logic
$(L, \mathcal{M})$, where $\mathcal{M}$ is a finite matrix, is given by:

1. all the sequents $\{A^i : i \in M\}$, for all formulas $A$, as axioms;

2. the weakening rule:
$$\frac{G}{G, F}$$

3. for each connective $c$ and each truth-value $i$, an introduction rule for $c$ at place $i$;

4. for any truth-values $i, j$, with $i \neq j$, the cut rule:
$$\frac{G, A^i \quad F, A^j}{G, F}$$

Of course, the above rules are sound, i.e., preserve satisfiability; for any structure $\sigma$, if $\sigma$ satisfies all the premises of a rule, then
it satisfies the conclusion. The axioms are valid because any formula must take at least one truth-value, and the cut rule is sound because no formula takes more than one value. Hence only valid sequents are provable. Conversely, every valid sequent is provable without cuts [3, Theorem 3.2].

These rules are far from unique. If we consider sequents as sequences of signed formulas, or multisets of signed formulas, then we need exchange and/or contraction rules. One might also allow introduction rules for $c$ at place $i$ in the form

$$\frac{G_1, F_1 \quad G_2, F_2 \quad \cdots \quad G_h, F_h}{G_1, G_2, \ldots, G_h, (cA_1 \ldots A_u)^i}$$

thus letting weakenings be implicit in the introduction rules.

A third possibility is to drop the weakening rules by taking all sequents of the form $G, \{A^i : i \in M\}$ as axioms, for $A$ any formula and $G$ any sequent. In this last case, the only rule in Definition 2.2.3 that remains not invertible is the cut rule (we recall that a rule is *invertible* if it preserves and reflects satisfiability: a structure satisfies all the premises iff it satisfies the conclusion). The cut elimination and the midsequent theorems are proved in [3].

Albeit not invertible, the following rules —together with axioms and weakening, but no cut— provide a complete sequent calculus for finite-valued logics [84], [98], [99]:

for any $u$-ary connective $c$, and any truth-values $i_1, \ldots, i_u$,

$$\frac{G, A_1^{i_1} \quad G, A_2^{i_2} \quad \cdots \quad G, A_u^{i_u}}{G, (cA_1 \ldots A_u)^i}$$

where $i = c(i_1, \ldots, i_u)$.

As a final reference, we cite the hypersequent approach to three-valued logics described in [1].

**Tableau calculi**

The sequent calculi introduced above do not need any elimination rule. For this reason they can be easily dualized, yielding semantic tableau calculi [96], [97], [13], [14]. As before, a sequent (in this context, we also say a *node*) is a finite set of signed formulas. Whereas in the sequent calculus we think of a sequent as the disjunction of its elements, in the tableau approach we think of it in conjunctive terms. We say that the structure $\sigma$ satisfies the node $G$ if $\overline{\sigma}(\bigwedge G) = 1$, i.e., for every $0 \leq i \leq n-1$ and every $A$ such that $A^i \in G$, we have $\sigma(A) = i$. $G$ is *unsatisfiable* if no $\sigma$ satisfies it.

Tableau calculi are refutational calculi; in order to prove $B_1, \ldots, B_r \models_\mathcal{M} A$, we prove the following:

for all mappings $\varphi : \{B_1, \ldots, B_r, A\} \to M$ such that $\varphi(B_1), \ldots, \varphi(B_r) \in D$ and $\varphi(A) \notin D$ (these maps are finite in number), the sequent $B_1^{\varphi(B_1)}, \ldots, B_r^{\varphi(B_r)}, A^{\varphi(A)}$ is unsatisfiable.

In particular, the formula $A$ is valid iff all the sequents $A^i$, for $i \in M \setminus D$, are unsatisfiable (in plain words, no structure assigns a nondesignated truth-value to $A$).

In order to prove that $G$ is unsatisfiable we construct a downward tree, whose root is $G$, and such that, for any node $F$ whose immediate successors are $F_1, \ldots, F_h$, and for any structure $\sigma$, $\sigma$ satisfies $F$ iff $\sigma$ satisfies at least one of $F_1, \ldots, F_h$. Rules for determining, given $F$, suitable immediate successors of $F$ can be formulated by dualizing the introduction rules for the connectives in the sequent calculus. Let $c$ be a connective of arity $u$, $i$ a truth-value, $C_1 \vee \cdots \vee C_h$ an $i$th DNF for $c$. Then $C_1 \vee \cdots \vee C_h$ induces the following *elimination rule for $c$ at place $i$*:

$$G, (cA_1 \ldots A_u)^i$$
$$\swarrow \quad \downarrow \quad \cdots \quad \searrow$$
$$G, F_1 \quad G, F_2 \quad \cdots \quad G, F_h$$

where $(cA_1 \ldots A_u)^i \notin G$ and $F_s$ is the sequent obtained from the conjunction $C_s$ (thought of as the set of its conjuncts) by substituting $A_r^j$ for $p_r^j$, for all $1 \leq r \leq u$ and all $j \in M$.

**Example 2.2.4** The 1/2th DNF for $\to$ in Example 2.2.1 induces the rule

$$G, (A \to B)^{1/2}$$
$$\swarrow \quad \searrow$$
$$G, A^{1/2}, B^0 \qquad G, A^1, B^{1/2}$$

For comparison's sake with the corresponding sequent calculus rule, we write the rule as

$$\Gamma \Rightarrow \Delta, A \to B \Rightarrow \Sigma$$

$$\Gamma, B \Rightarrow \Delta, A \Rightarrow \Sigma \qquad \Gamma \Rightarrow \Delta, B \Rightarrow \Sigma, A$$

It is clear that, whichever sequent we take as the root, the process of constructing the tree cannot run indefinitely, and we will reach a position where no rule can be further applied. This happens when every leaf $F$ in the tree is *terminal*, i.e., the following hold:

- either the only signed formulas in $F$ are literals;

- or every signed formula in $F$ which is not a literal is of the form $(cA_1 \ldots A_u)^i$, for a connective $c$ and a truth-value $i$ for which there is no elimination rule for $c$ at place $i$ (i.e., the empty disjunct is an $i$th DNF for $c$ or, equivalently, $c$ never takes the value $i$).

Given a terminal leaf $F$, we say that the corresponding branch is *open* if $F$ contains only literals and, for no $p, i, j$ with $i \neq j$, $F$ contains both $p^j$ and $p^i$; otherwise, we say that the branch is *closed*. Assume we have a tree whose root is $G$ and whose leaves are all terminal. Then, for any $\sigma$, $\sigma$ satisfies $G$ iff $\sigma$ satisfies at least one of the leaves. $G$ is unsatisfiable iff all branches are closed, whereas, if $F$ is a leaf whose branch is open, then any $\sigma$ such that $(\sigma(p) = i$ iff $p^i \in F)$ satisfies $G$.

These ideas are trivially modifiable in various ways; it is usual to write the elimination rule for $c$ at place $i$ in the form

$$G, (cA_1 / \ldots A_u)^i$$

$$F_1 \quad F_2 \qquad F_h$$

where the slash sign means that the occurrence of $(cA_1 \ldots A_u)^i$ has been *used* and does not appear any more in the branches through it. We close a branch as soon as a formula appears in it signed by two different truth-values, and we take as the final leaf the set of all signed formulas still appearing along the branch.

In [39], Hähnle introduces the idea of using sets of truth-values as signs. Define a *sign* to be a nonempty subset $S$ of $M$. A *signed formula* is now an expression of the form $A^S$; its intended meaning is that the formula $A$ takes a truth-value in $S$. We say that the structure $\sigma$ satisfies $A^S$ if $\sigma(A) \in S$; the old signed formula $A^i$ correspond to the new $A^{\{i\}}$.

Let $\mathcal{P}(M)$ be the algebra, of the same similarity type as $M$, whose elements are the nonempty subsets of $M$. If $c$ is a connective of arity $u$ and $S_1, \ldots, S_u \in \mathcal{P}(M)$, then $c(S_1, \ldots, S_u)$ is defined to be $\{c(i_1, \ldots, i_u) : i_t \in S_t\}$. Choose a subalgebra $\mathcal{S}$ of $\mathcal{P}(M)$; $\mathcal{S}$ will constitute the set of signs of our tableau calculus. For any $S \in \mathcal{S}$ and any connective $c$, one constructs an $S$th DNF for $c$, meaning a signed formula expression of the form $C_1 \vee \cdots \vee C_h$, where each $C_s$ is a conjunct of signed formulas from the set $\{p_r^T : 1 \leq r \leq u, T \in \mathcal{S}\}$. The condition to be fulfilled is that, for any $\sigma$, $\sigma$ satisfies $(cp_1 \ldots p_u)^S$ iff there exists $s$ such that $\sigma$ satisfies all the signed formulas in $C_s$. As before, $C_1 \vee \cdots \vee C_h$ induces a tableau elimination rule for $c$ at place $S$.

**Example 2.2.5** Consider Gödel's three-valued logic of §2.5.2. We want to construct an elimination rule for $\wedge$ at place $1/2$. In the old formulation we have the rule

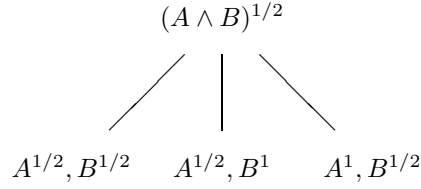

Let $\mathcal{S}$ be the subalgebra of $\mathcal{P}(\mathcal{G}_3)$ generated by $\{\{1/2\}, \{1/2, 1\}\}$. Then we obtain the rule

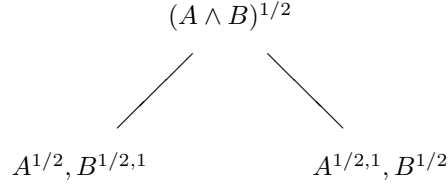

The construction of a tableau gets modified in a straightforward fashion. For example, we close a branch as soon as:

- either a formula $(cA_1 \ldots A_u)^S$ appears along the branch and no elimination rule for $c$ at place $S$ is defined;

- or the branch contains signed formulas $A^{S_1}, \ldots, A^{S_h}$ with $S_1 \cap \cdots \cap S_h = \emptyset$.

The resulting system is sound and complete [40]. In order to prove that the formula $A$ is valid, we do not need any more to construct, for each nondesignated truth-value $i$, a closed tableau for $A^i$, but a single closed tableau for $A^{M \setminus D}$ will do (provided $M \setminus D \in \mathcal{S}$). What are the advantages of the sets-of-signs method? The new rule in Example 2.2.5 not only has fewer conclusions, but it is a $\beta$-rule, in the sense of [89]. If for a certain class of logics we can formulate tableau systems using only $\alpha$ and $\beta$ rules, then the resulting uniform notation makes it possible to prove certain desired properties in a uniform manner. The disadvantage of the method is of course the increase in the number of signs. Hähnle introduces the class of *regular connectives* (roughly, connectives which —in the sets-as-signs presentation— admit an $\alpha$ or $\beta$ tableau elimination rule), and proves that the class of finite-valued logics expressible in terms of regular connectives only is quite large. As a matter of fact, for each $n$, the functionally complete Post $n$-valued logic has a formulation in terms of regular connectives [40, p. 69]. See also [38] for an up-to-date introduction to tableau methods in many-valued logics.

### 2.2.2 Resolution calculi

Let $\mathcal{M} = (M, D)$ be a fixed finite matrix. Suppose that it is possible to construct, for every finite set of formulas $\Gamma$ and every formula $A$, a new finite set of formulas $R(\Gamma, A)$ such that

$$\Gamma \models_{\mathcal{M}} A \quad \text{iff} \quad R(\Gamma, A) \text{ is unsatisfiable.}$$

Then the problem of determining the validity of $\Gamma \models_{\mathcal{M}} A$ is reduced to the problem of the satisfiability of $R(\Gamma, A)$. Observe that, since we are dealing with finite matrices, the relation $\Gamma \models_{\mathcal{M}} A$ is decidable. Hence, it is always possible to construct an appropriate $R(\Gamma, A)$, and the above reduction makes sense only if passing from $(\Gamma, A)$ to $R(\Gamma, A)$ is "easy"; this is the case for most common logics. A sufficient condition is the following:

- there exists a formula $N(p) \in FORM(L)$ whose only sentential variable is $p$, and such that, for every $i \in M$,

$$i \in D \quad \text{iff} \quad N(i) \notin D.$$

Indeed, assume that such a formula exists. Then we have

$$\begin{aligned}
B_1,\ldots,B_r \models_{\mathcal{M}} A \quad &\text{iff} \quad \forall \sigma(\sigma(\{B_1,\ldots,B_r\}) \not\subseteq D \text{ or } \sigma(A) \in D) \\
&\text{iff} \quad \forall \sigma(\sigma(\{B_1,\ldots,B_r\}) \not\subseteq D \text{ or } \sigma(N(A)) \notin D) \\
&\text{iff} \quad \{B_1,\ldots,B_r, N(A)\} \text{ is unsatisfiable.}
\end{aligned}$$

**Example 2.2.6** For classical logic, just let $N(p) = \neg p$. Functional completeness of Post's $n$-valued logic (see §2.3) guarantees the existence of an appropriate $N(p)$, while for Łukasiewicz's $n$-valued logic one can take

$$N(p) = p \to (p \to (p \to \cdots \to (p \to \neg p)))$$

($n-1$ implication signs).

Let $\Gamma$ be a finite set of formulas; we want to use the resolution rule to determine if $\Gamma$ is satisfiable. There are two main approaches: the clausal one and the non-clausal.

**The clausal approach**

Let $A$ be a formula, and let $\{i_1,\ldots,i_r\}$ be the set of designated truth-values. We construct a signed formula expression $E_A$ such that:

- $E_A$ is in CNF, does not contain the connective $\neg$, and contains only signed formulas of the form $p^j$, for $p$ a propositional variable appearing in $A$, and $j$ a truth-value;

- $A^{i_1} \vee \cdots \vee A^{i_r} \equiv E_A$.

Such a signed formula expression can be constructed by induction on the complexity of $A$, using appropriate $i$th conjunctive and disjunctive normal forms for the connectives. Boolean equivalences can be used freely, and moreover:

- subexpressions of the form $B^i \wedge B^j$, with $i \neq j$, may be replaced by *false*;

- subexpressions of the form $B^{i_1} \vee \cdots \vee B^{i_n}$, with $\{i_1,\ldots,i_n\} = $ the set of all the truth-values, may be replaced by *true*.

**Example 2.2.7** Let us work in Łukasiewicz's three-valued logic. Then:

$$\begin{aligned}
(p \to (p \to \neg p))^1 &\equiv p^0 \vee (p \to \neg p)^1 \vee (p^{1/2} \wedge (p \to \neg p)^{1/2}) \\
&\equiv p^0 \vee p^0 \vee (\neg p)^1 \vee (p^{1/2} \wedge (\neg p)^{1/2}) \vee (p^{1/2} \wedge (p^{1/2} \vee (\neg p)^{1/2}) \wedge (p^0 \vee (\neg p)^0)) \\
&\equiv p^0 \vee p^{1/2} \vee (p^{1/2} \wedge (p^0 \vee p^1)) \\
&\equiv (p^0 \vee p^{1/2}) \wedge (p^0 \vee p^{1/2} \vee p^1) \\
&\equiv p^0 \vee p^{1/2}
\end{aligned}$$

Let $\Gamma$ be a finite set of formulas. In order to test the satisfiability of $\Gamma$ we proceed in two steps:

1. for each $A \in \Gamma$, we construct $E_A$ as above, and we set $E_\Gamma$ to be the set of all disjunctions that appear as conjuncts in some $E_A$. As it is usual in the resolution literature, we drop the connective $\vee$, and we regard $E_\Gamma$ as a set of clauses. $E_\Gamma$ can be seen as a logic-free equivalent —with respect to satisfiability— to $\Gamma$;

2. we apply the resolution rule to the clauses in $E_\Gamma$. The rule has the form

$$\frac{D, p^i \quad D', p^j}{D, D'}$$

where $D, p^i$ and $D', p^j$ are previously deduced clauses, and $i \neq j$.

It is not difficult to prove the completeness of the resolution rule: a finite set of clauses is unsatisfiable iff the empty clause is derivable from it. As the set of clauses derivable from a finite set is finite (no new literal is ever added), this yields a decision procedure for satisfiability. We refer to [62], [2], [71], [41] for further references and the extension to the predicate case.

**The non-clausal approach**

In non-clausal resolution systems we do not convert the set $\Gamma$ into clausal form. Instead, we consider a pair $(V, \mathcal{F})$ such that:

- $V$ is a finite algebra, of the same similarity type as the algebra $M$;
- the elements of $V$ —which are called *verifiers*— are formulas in the language $L$;
- $\mathcal{F}$ is a family of subsets of $V$, and a subset of $V$ is unsatisfiable iff it extends an element of $\mathcal{F}$.

We fix an enumeration $W_1, \ldots, W_h$ of the elements of $V$. The resolution rule induces a branching process, and it is applied to $\Gamma$ in the form

$$A_1(p), A_2(p), \ldots, A_h(p)$$
$$A_1(W_1) \quad A_2(W_2) \quad A_h(W_h)$$

where $A_1(p), \ldots, A_h(p)$ are previously deduced formulas, not necessarily distinct, sharing a common variable $p$ which does not occur in $W_1, \ldots, W_h$. The idea behind the rule is that $\{A_1(p), \ldots, A_h(p)\}$ is satisfiable iff, for some $1 \leq i \leq h$, $\{A_1(p), \ldots, A_h(p), A_s(W_s)\}$ is satisfiable.

In addition to the resolution rule, we have *transformation* rules. These are simply the tables defining the connectives in $V$, so that if $c$ is a connective of arity $u$, $T_1, \ldots, T_u \in V$, and $c(T_1, \ldots, T_u) = T$ in $V$, then we can simplify $A(cT_1, \ldots, T_u)$ to $A(T)$.

The computation proceeds by constructing a tree, whose nodes are labelled by finite sets of formulas. The root of the tree is labelled by $\Gamma$, and branching is induced by the resolution rule. The transformation rules do not induce branching, and the computation reports that $\Gamma$ is unsatisfiable if a stage is reached where, for every branch of the tree, the set of formulas occurring along the branch extends to an element of $\mathcal{F}$. Hence, the rôle of the family $\mathcal{F}$ is to provide the termination conditions for the deductive process.

**Example 2.2.8** In classical logic, as formulated in Example 2.1.1, we set $V = \{F, T\}$, where $F$ is any contradictory formula and $T$ any tautology, and $\mathcal{F} = \{\{F\}\}$. If we want to show the unsatisfiability of $\Gamma = \{p, q, \neg p \vee \neg q\}$, we construct the tree

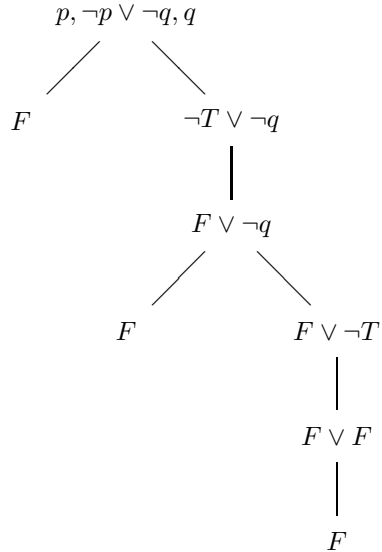

When the language contains a disjunction connective, then the resolution rule can take the form of a single conclusion rule, and the deductive process does not require branching [92].

A positive aspect of Stachniak's approach is its applicability to a class of consequence operations larger than the class of operations definable by strongly adequate matrices. In particular, it works for operations $C$ for which there exists a finite set $\mathcal{K}$ of finite matrices such that, for every $\Gamma$,

$$C(\Gamma) = \bigcap \{C_\mathcal{M}(\Gamma) : \mathcal{M} \in \mathcal{K}\}.$$

Such a consequence operation is called *strongly finite*. The negative aspect is that the cardinality of $V$ will be in general much greater than the cardinality of the matrix $\mathcal{M}$ (assuming that we are working with a single-matrix consequence). For example, six verifiers are needed for Łukasiewicz's three-valued logic. See [95], [93], [92] for further reference; see also the discussion in [40, §8.1.3].

## 2.3 Post systems

Post introduced $n$-valued systems in his doctoral thesis, which was published in [78]. His approach is purely algebraic: he develops his calculi as a direct generalization of classical 2-valued calculus. He gives no linguistic interpretation to his truth-values, except in terms of sequences of ordinary truth-values, anticipating the monotonic representation of elements in a Post algebra (see below). After his 1921 paper, he never came back to many-valued logic.

**Definition 2.3.1** Let $n \geq 2$, and let $0 \leq m \leq n-1$. The matrix $\mathcal{P}_n^m$ for $n$-valued Post logic is defined as follows: $\mathcal{P}_n^m = (\mathbf{n}, \{m, m+1, \ldots, n-1\})$. There are only two connectives, $\vee$ and $\sim$, defined by $i \vee j = \max\{i, j\}$, and $\sim i = i - 1 \pmod{n}$. When we are not concerned with the set of designated values, we denote the algebra $(\mathbf{n}, \vee, \sim)$ by $P_n$.

Classical negation in $\mathbf{2}$ can be seen either as mirror reflection ($i \mapsto 1 - i$), or as cyclic shift ($i \mapsto i - 1 \pmod{2}$); the two things coincide. For $n \geq 3$ they do not coincide any more; Post's negation $\sim$ takes the second alternative.

The key feature of $P_n$ is functional completeness: for any $k \geq 1$, any function $f : \mathbf{n}^k \to \mathbf{n}$ is definable in terms of $\vee$ and $\sim$ (in particular $i \wedge j = \min\{i, j\}$ is definable). This was proved by Post himself; another proof can be found in [104]. We want to introduce all 0-ary connectives. Of course, it suffices to introduce the constant 0. We may define 0 to be $p \wedge \sim p \wedge \sim\sim p \wedge \cdots \wedge \sim^{n-1} p$, or we may introduce 0 as a new primitive connective. Either way, we are able to allow $k$ in the definition of functional completeness to be 0 as well. As —abstractly— it does not matter which set of connectives we assume as primitive, as long as they define the same class of functions, we are justified in calling $n$-valued Post logic any $n$-valued logic which is functionally complete ($k = 0$ allowed). We note incidentally that from the practical point of view —say, from the point of view of automated theorem proving— it makes a lot of difference which set of connectives we assume as primitive. A bad choice of primitives can make work in a system unpractical: nobody has ever tried to set up classical logic by using the Sheffer stroke only.

A remarkable feature of $n$-valued Post logic is Craig's interpolation theorem.

**Theorem 2.3.2** *Assume that the sequent $A_0 \Rightarrow A_1 \Rightarrow \cdots \Rightarrow A_{n-1}$ is valid. Then there is a formula $B$ such that:*

1. *any propositional symbol occurring in $B$ occurs in at least two of the $A_i$'s;*

2. *all the sequents $A_0 \Rightarrow B \Rightarrow \cdots \Rightarrow B$, $B \Rightarrow A_1 \Rightarrow \cdots \Rightarrow B$, ..., $B \Rightarrow B \Rightarrow \cdots \Rightarrow A_{n-1}$ are valid.* ∎

The above theorem was proved for first-order Post logic by Hanazawa and Takano in [46], generalizing results in [32] and [61]; see also [106], [72], [114].

On the algebraic side, much effort has been devoted to the study of the class of Post algebras. The main ideas are as follows: we give $P_n$ the similarity type $(2, 2, 1, 0, 0)$, in the language $(\vee, \wedge, \sim, \bot, \top)$, with $\bot$ and $\top$ defined to be 0 and $n-1$, respectively. As an abstract algebra, $P_n$ generates an equational class, which we denote by $\mathbf{P}_n^\sim$. The elements of $\mathbf{P}_n^\sim$ are all the algebras, of the same similarity type of $P_n$, that satisfy the identities satisfied by $P_n$. By Birkhoff's theorem, $\mathbf{P}_n^\sim$ is the class of all homomorphic images of all subalgebras of all products of copies of $P_n$; $\mathbf{P}_2^\sim$ is the class of Boolean algebras.

We now drop $\sim$. By definition, the class $\mathbf{P}_n$ of *Post algebras of order $n$* is the class of algebras whose elements are the elements of $\mathbf{P}_n^\sim$, but in the language $(\vee, \wedge, \bot, \top)$. $\mathbf{P}_n$ is no longer an equational class, since it is not closed for subalgebras. An algebra $N \in \mathbf{P}_n$ is a distributive lattice with bottom and top, and with an extra feature: $N$ has buried in its structure the chain of constants $c_{n-1} = \top, c_{n-2} = \sim$ $\top, c_{n-3} =\sim\sim \top,\ldots, c_0 =\sim^{n-1} \top = \bot$, even though it does not remember them any more. Can we recover this chain?

We recall that an element $a$ of a distributive lattice with bottom and top is *complemented* if there exists $b$ such that $a \vee b = \top$ and $a \wedge b = \bot$. If such a $b$ exists, then it is uniquely determined; we denote it by $\neg a$. By definition, the *center* $\mathcal{C}(N)$ of the lattice $N$ is the Boolean algebra of all complemented elements of $N$.

**Theorem 2.3.3** [26] *Let $N = (N, \vee, \wedge, \bot, \top)$ be a distributive lattice with bottom and top. Then $N$ is a Post algebra of order $n$ if and only if it contains a chain $\bot = c_0 < c_1 < c_2 < \cdots < c_{n-1} = \top$ such that:*

1. *every element $a \in N$ can be written as*
$$a = (a_1 \wedge c_1) \vee (a_2 \wedge c_2) \vee \cdots \vee (a_{n-1} \wedge c_{n-1}),$$
*with $a_1, \ldots, a_{n-1} \in \mathcal{C}(N)$;*

2. *for no $a \in \mathcal{C}(N)$ except $\bot$, and no $0 \leq i \leq n-2$, it is $a \wedge c_{i+1} \leq c_i$.* ∎

If we require that $a_1 \geq a_2 \geq \cdots \geq a_{n-1}$ then the $a_i$'s are uniquely determined; they constitute the *monotonic representation* of $a$. The chain of the $c_i$'s, as well as its length, is uniquely determined; it follows that, if $n \neq m$, the intersection of $\mathbf{P}_n$ and $\mathbf{P}_m$ is empty.

The situation become clearer as soon as one knows about Epstein's version of the Stone representation theorem.

**Theorem 2.3.4** [27] *Any Post algebra of order $n$ is isomorphic to the lattice of all continuous functions from the dual space of its center —a compact, totally disconnected, Hausdorff space— to the $n$-element chain. Conversely, for any Boolean space $X$, the lattice of all continuous functions : $X \to \mathbf{n}$ is a Post algebra of order $n$.* ∎

By Speed's theorem [90], this amounts to saying that every Post algebra of order $n$ is the coproduct of a Boolean algebra —namely, its center— with the $n$-element chain.

In Epstein's representation, the constant $c_i$ corresponds to the function taking constantly the value $i$, whereas $a_i$ in the monotonic representation of $a$ is the characteristic function of the set $\{x \in X : a(x) \geq i\}$. It is also clear that the monotonic representation yields an isomorphism of $N$ onto a substructure $Q$ of the product of $n-1$ copies of $\mathcal{C}(N)$. $Q$ is formed by taking all decreasing sequences of length $n-1$ of elements of $\mathcal{C}(N)$. The operations $\vee, \wedge, \top, \bot$ act componentwise, while the shift $\sim$ is given by:
$$(a_1, a_2, \ldots, a_{n-2}, a_{n-1}) \mapsto (a_2 \vee \neg a_1, a_3 \vee \neg a_1, \ldots, a_{n-1} \vee \neg a_1, \neg a_1)$$

It is remarkable that at the end of his 1921 paper, Post anticipated this last representation, by interpreting his truth-values as monotonic sequences of ordinary truth-values.

## 2.4 Łukasiewicz systems

As we saw in §2.3, Post systems are quite rigid: in order to pass from $n$ to $m$ values, one has to change the chain of constants. Łukasiewicz systems are more flexible; we already encountered the three-valued version in Example 2.1.2, and we now give the general framework.

Let $n \geq 2$. The matrix $\mathcal{L}_n$ for $n$-valued Łukasiewicz logic is $\mathcal{L}_n = (I_n, \{1\})$, in the language $(\to, \neg, 1)$. The constant 1 is to be interpreted in the real number 1, while $\to$ and $\neg$ are defined thus: $i \to j = \min\{1 - i + j, 1\}$, $\neg i = 1 - i$. If we substitute $I_n$ with the real unit interval $I = [0, 1]$, with no changes in the definitions of the connectives, we obtain the matrix $\mathcal{L} = (I, \{1\})$ for Łukasiewicz infinite-valued logic.

We introduce the following connectives:
$$\begin{aligned}
0 &= \neg 1 \\
i \oplus j &= (\neg i) \to j \\
i \otimes j &= \neg((\neg i) \oplus (\neg j)) \\
i \vee j &= (i \to j) \to j \\
i \wedge j &= \neg((\neg i) \vee (\neg j)) \\
i \leftrightarrow j &= (i \to j) \wedge (j \to i)
\end{aligned}$$

We easily see that:

- $\oplus, \otimes, \vee, \wedge$ are all associative and commutative;
- $\vee$ and $\wedge$ correspond to max and min;
- $i_1 \oplus \cdots \oplus i_t = \min\{i_1 + \cdots + i_t, 1\}$;
- $i_1 \otimes \cdots \otimes i_t = \max\{i_1 + \cdots + i_t - t + 1, 0\}$;
- $i \leftrightarrow j = (i \to j) \otimes (j \to i) = 1 - |i - j|$;
- if restricted to $\{0, 1\}$ the connectives $\to, \neg, \vee, \wedge, \leftrightarrow$ behave in a Boolean way, while $\oplus$ and $\otimes$ coincide with $\vee$ and $\wedge$, respectively;
- the adjunction $(i \leq j \to k$ iff $i \otimes j \leq k)$ holds;
- the sets of connectives $(\to, \neg, 1)$, $(\to, 0)$, $(\oplus, \otimes, \neg, 0, 1)$ are interdefinable.

In the literature the weight of Łukasiewicz logic has always been predominant. Undeniably, tradition played a rôle: Łukasiewicz logic was the first many-valued logic to be proposed and systematically studied [55], [56]. More recently, new motivations arose with the development of fuzzy logic. Indeed, it seems reasonable that a logic aiming at the formalization of any inferential mechanism should come equipped with a "product" connective $\otimes$ and an "implication" one $\to$, related by the adjunction condition. Following [22], [77], we define a *residuated lattice* to be a structure $(L, \vee, \wedge, \otimes, \to, \bot, \top)$ such that:

- $(L, \vee, \wedge, \bot, \top)$ is a lattice with bottom and top;
- $(L, \otimes, \top)$ is a commutative monoid;
- the adjunction condition holds, i.e., for any $i, j, k \in L$, we have $i \leq j \to k$ iff $i \otimes j \leq k$.

If, in addition, $L$ is complete as a lattice, we say that $(L, \vee, \wedge, \otimes, \to, \bot, \top)$ is a *complete residuated lattice*.

**Theorem 2.4.1** [77, pp. 121-122] *Let $(L, \vee, \wedge, \bot, \top)$ be a complete lattice.*

1. *Let $\otimes$ be a binary operation (not necessarily commutative) over $L$ such that:*

   a. *$\otimes$ is isotone in both variables (i.e., $i \leq j$ implies $i \otimes k \leq j \otimes k$ and $k \otimes i \leq k \otimes j$);*

   b. *for each subset $I$ of $L$, $(\bigvee_{i \in I} i) \otimes j = \bigvee_{i \in I} (i \otimes j)$.*

   *Then there exists a unique operation $\to$ over $L$ such that the adjunction condition holds for the pair $(\otimes, \to)$; this operation is given by $j \to k = \bigvee\{i : i \otimes j \leq k\}$.*

2. *Let $\to$ be a binary operation over $L$ such that:*

   c. *$\to$ is antitone in the first variable and isotone in the second;*

   d. *for each subset $K$ of $L$, $j \to (\bigwedge_{k \in K} k) = \bigwedge_{k \in K} (j \to k)$.*

   *Then there exists a unique operation $\otimes$ over $L$ such that the adjunction condition holds for the pair $(\otimes, \to)$; this operation is given by $i \otimes j = \bigwedge\{k : i \leq j \to k\}$.*

3. *In every complete residuated lattice the conditions (a), (b), (c), (d) hold.* ∎

**Example 2.4.2** Consider the ordinary lattice structure of the real unit interval $I = [0, 1]$ As proved in [77, p. 123] and references therein, $I$ can be given $2^{\aleph_0}$ distinct structures of residuated lattice in which the product operation is continuous.

1. Take $\otimes = \wedge$. Then $\to$ is defined by:

$$i \to j = \begin{cases} 1, & \text{if } i \leq j; \\ j, & \text{otherwise.} \end{cases}$$

   We get a complete Heyting algebra. This is the case of Gödel's logic; see §2.5.2.

2. Let $\otimes$ be the ordinary product between real numbers. Then we get:

$$i \to j = \begin{cases} 1, & \text{if } i \leq j; \\ j/i, & \text{otherwise.} \end{cases}$$

See [44] for this product logic.

3. Łukasiewicz's product $i \otimes j = (i+j-1) \vee 0$ yields Łukasiewicz's implication $i \to j = (1-i+j) \wedge 1$.

These three examples are in some sense exhaustive; indeed, every continuous $t$-norm (see Novák's Chapter of this Handbook) can be obtained from them [76].

**Theorem 2.4.3** [60] *Let $(I, \vee, \wedge, \otimes, \to, 0, 1)$ be a residuated lattice in which $\vee$ and $\wedge$ are the usual* max *and* min. *Suppose that $\to: I^2 \to I$ is continuous. Then $\otimes$ and $\to$ are the Łukasiewicz product and implication.* ∎

The above theorem is even more surprising if one considers that 56 years separate [55] and [60]. Alone, it would be sufficient to justify the study of Łukasiewicz logic. As a matter of fact, Pavelka proved in [77] that the Łukasiewicz connectives guarantee a sort of uniqueness with respect to the completeness theorem in a fuzzy logic context.

These considerations are relatively recent; Łukasiewicz logic has not been born as the logic of residuated lattices, but as the logic of a certain matrix semantics. How can we axiomatize the set $C_{\mathcal{L}}(\emptyset)$ of all formulas that are valid in $\mathcal{L}$? Łukasiewicz conjectured that the following four axiom schemas, along with modus ponens, would do:

Ax1. $\alpha \to (\beta \to \alpha)$;

Ax2. $(\alpha \to \beta) \to ((\beta \to \gamma) \to (\alpha \to \gamma))$;

Ax3. $(\neg \alpha \to \neg \beta) \to (\beta \to \alpha)$;

Ax4. $((\alpha \to \beta) \to \beta) \to ((\beta \to \alpha) \to \alpha)$.

(actually, he also had $(\alpha \to \beta) \vee (\beta \to \alpha)$ as a fifth axiom; this was proved to be redundant by Chang and Meredith). In 1935 Wajsberg claimed to have a proof of the conjecture, but he never published it. In 1958 Rose and Rosser [82] gave the first published proof; other proofs have been obtained by Chang [16], Cignoli [18], and Panti [75]. We will return to the completeness issue in §3.

Axioms Ax1–Ax4 are quite elegant. If we add to them

Ax5. $(\neg \alpha \to \alpha) \to \alpha$;

we obtain a variant of the Hilbert-Ackermann axiomatization of the two-valued calculus. If we add

Ax5'. $((\alpha \to \neg \alpha) \to \alpha) \to \alpha$;

we get an axiomatization of Łukasiewicz's three-valued calculus; in both cases Ax4 becomes redundant and may be dropped.

Ax5 expresses the idempotence of $\oplus$; in fact, it can be equivalently formulated as $\alpha \oplus \alpha \to \alpha$ ($\alpha \to \alpha \oplus \alpha$ being deducible from Ax1–Ax4). Under Ax1–Ax5, $\oplus$ and $\vee$ collapse, and so do $\otimes$ and $\wedge$. For $n \geq 3$, the axiom parallel to Ax5 is

Ax5''. $\underbrace{\alpha \oplus \alpha \oplus \cdots \oplus \alpha}_{n \text{ summands}} \to \underbrace{\alpha \oplus \cdots \oplus \alpha}_{n-1 \text{ summands}}$.

For $n = 3$, this is enough (i.e., Ax1–Ax3,Ax5' are equivalent to Ax1–Ax4,Ax5'' and axiomatize $C_{\mathcal{L}_3}(\emptyset)$), but for $n \geq 4$ we further need an entire group of new axioms, namely

Ax6j. $(n-1)(\underbrace{\alpha \oplus \alpha \oplus \cdots \oplus \alpha}_{j \text{ summands}} \to (\alpha \otimes (\underbrace{\alpha \oplus \cdots \oplus \alpha}_{j-1 \text{ summands}})))$;

(see [37]). Here $(n-1)\beta$ stands for $\beta \oplus \cdots \oplus \beta$ ($n-1$ summands), while $j$ varies in the set $\{k : 1 < k < n-1$ and $k$ does not divide $n-1\}$.

What remains open is the problem of a "linguistic" interpretation of the logic. The most intriguing connective is the truncated sum $\oplus$. To compute $i \oplus j$, add $i$ and $j$ as real numbers: if the result is less or equal to 1, take it, otherwise, take 1. Are there situations in real life in which we act on proposition by truncated addition? In [87], Scott suggests that *degrees of error* may be considered additive. So if $A$ is true to degree of error $i$ and $B$ is true to degree of error $j$ ($i, j \in [0, 1]$), then ($A$ and $B$) should be true to degree of error $\min\{i + j, 1\}$ (degrees of error greater than 1, i.e., complete falsity, can of course be cut down to 1). But "$A$ is true to degree of error $i$" in the structure $\sigma$ means that $\sigma(A) = 1 - i$. Hence we get:

$$\begin{aligned}
\sigma(A \text{ and } B) &= 1 - ((i + j) \wedge 1) \\
&= 1 + ((-i - j) \vee -1) \\
&= (1 - i - j) \vee 0 \\
&= (1 - i + 1 - j - 1) \vee 0 \\
&= (\sigma(A) + \sigma(B) - 1) \vee 0 \\
&= \sigma(A) \otimes \sigma(B)
\end{aligned}$$

which is indeed what we wanted; the adjunction condition now gives us the definition of $\rightarrow$.

As a by-product, we have the non-idempotence of conjunction: we must be willing to accept that, if $A$ is true to degree of error, say, $2/3$, then ($A$ and $A$) is true to degree of error $(2/3) \otimes (2/3) = 1/3$. We will return to Scott's interpretation in the final part of this survey; in the meantime, we refer to Scott's and Smiley's discussion in [87], as well as to Urquhart's [104]. We also refer to [87] for Giles's interpretation of Łukasiewicz logic as a "logic of risk and commitment". We only remark that the non-idempotence of "and" is not so shocking: one might think of two different trials of the same experiment, or of two different statements about the same trial. This is not so far from, e.g., linear logic.

In [67], Mundici proposes an interpretation of Łukasiewicz logic in terms of the Ulam game with lies. The situation is as follows: an Answerer chooses a number between 0 and $k$. A Questioner asks yes-no questions. The Answerer is allowed to lie up to $n$ times. How many questions does the Questioner need to find the number? Neither the principle of non-contradiction nor idempotence hold in this context. If the Answerer says first "The number is 7", and then "The number is not 7", this does not lead to inconsistency: it just means that the Answerer has one lie less at his disposal. Similarly, repeated assertions that "The number is 7" are more informative that one single assertion; $n + 1$ assertions guarantee truth. Mundici shows that, at each stage of the game, both the Questioner state of knowledge and the Answerer's assertion can be expressed by a formula in the $k + 2$-valued Łukasiewicz logic. The Questioner's next state of knowledge is then given by the Łukasiewicz conjunction of the two formulas. The situation can be generalized by allowing an unbounded number of lies, or random lies. The Questioner's adaptive strategy leads naturally to the theory of error-correcting codes; we refer to [67], [68] for further references.

There are two situations in which truncated addition arises naturally; I am grateful to S. Stefani for first calling my attention to these. The first situation is saturation in a logical gate. There exist logical gates which return as output the sum of the voltages they receive as input; the voltages themselves can be seen as truth-values. As every real-life electronic component, these gates have a saturation level, above which either the output does not increase any more, or the gate shorts out.

The second situation if fixed-point addition in computers. Fixed-point addition is addition $\mod 2^{32}$, or $\mod 2^{64}$, depending on the computer. It is less resource-consuming than floating-point addition, but it is subject to overflow. As the overflow generally results in a disaster, one may try to limit the damages by always chopping the result of a potentially dangerous addition to the maximum number expressible ($2^{32} - 1$, or $2^{64} - 1$).

## 2.5 Other systems

### 2.5.1 Bochvar's and Kleene's systems

Bochvar's and Kleene's are both three-valued systems. In addition to 0 and 1 for *false* and *true*, they have a third value 2. While for Łukasiewicz the third value stands for *possible*, or *not yet determined*, from Bochvar's point of view it stands for *paradoxical*, or *meaningless*. Any compound proposition that includes a meaningless part is meaningless itself, and hence the truth-tables

|   | $\neg$ |
|---|---|
| 0 | 1 |
| 1 | 0 |
| 2 | 2 |

| $\vee$ | 0 | 1 | 2 |
|---|---|---|---|
| 0 | 0 | 1 | 2 |
| 1 | 1 | 1 | 2 |
| 2 | 2 | 2 | 2 |

| ∧ | 0 | 1 | 2 |
|---|---|---|---|
| 0 | 0 | 0 | 2 |
| 1 | 0 | 1 | 2 |
| 2 | 2 | 2 | 2 |

| → | 0 | 1 | 2 |
|---|---|---|---|
| 0 | 1 | 1 | 2 |
| 1 | 0 | 1 | 2 |
| 2 | 2 | 2 | 2 |

Bochvar's systems was proposed in [11] as a way for avoiding the logical paradoxes, notably Russell's paradox. We refer to [81, §2.4] and [104, §1.6] for a deeper analysis and further references.

In [51], [52], Kleene regards the value 2 as *indeterminate*, in the same sense as a Turing machine which does not halt on a specific input can be seen as giving an indeterminate output. Two points of view are possible. First possibility: one assumes that an indeterminate partial output makes any computation globally indeterminate. With this understanding, one obtains exactly Bochvar's truth-tables: Kleene calls "weak connectives" the connectives defined in this way. Second possibility: one accepts that a computation may assume a determinate global value, even if some subcomputations are globally indeterminate. In this case, Kleene proposes what he calls the "strong connectives".

| | ¬ |
|---|---|
| 0 | 1 |
| 1 | 0 |
| 2 | 2 |

| ∨ | 0 | 1 | 2 |
|---|---|---|---|
| 0 | 0 | 1 | 2 |
| 1 | 1 | 1 | 1 |
| 2 | 2 | 1 | 2 |

| ∧ | 0 | 1 | 2 |
|---|---|---|---|
| 0 | 0 | 0 | 0 |
| 1 | 0 | 1 | 2 |
| 2 | 0 | 2 | 2 |

| → | 0 | 1 | 2 |
|---|---|---|---|
| 0 | 1 | 1 | 1 |
| 1 | 0 | 1 | 2 |
| 2 | 2 | 1 | 2 |

The table for implication does not follow from the one for conjunction, plus an adjunction condition, but is explicitly stated. Kleene defines classically $i \to j$ as $\neg i \vee j$, i.e., he considers the computation $i \to j$ successful if either $j$ is successful or $i$ unsuccessful, unsuccessful if $i$ is successful and $j$ unsuccessful, and indeterminate in all other cases. See [81], [104], as well as Kleene's original papers; see also [5] for the basic theory of Kleene algebras.

### 2.5.2 Gödel's system

Let $\mathcal{G}_n$ be the matrix $\mathcal{G}_n = (I_n, \{1\})$ in the language $(\vee, \wedge, \to, \neg, 0, 1)$. Here $\vee, \wedge$ are the max and min functions for the usual ordering of $I_n$. The implication $\to$ is given by adjunction with respect to $\wedge$ as in Example 2.4.2, while $\neg i$ is $i \to 0$. Explicitly,

$$i \to j = \begin{cases} 1, & \text{if } i \leq j; \\ j, & \text{otherwise}. \end{cases} \qquad \neg i = \begin{cases} 1, & \text{if } i = 0; \\ 0, & \text{otherwise}. \end{cases}$$

$\mathcal{G}_n$ may be extended to an infinite-valued matrix by replacing $I_n$ with any infinite subset of the real unit interval $[0, 1]$, with no change in the definitions of the connectives. Any two such extensions induce the same set of valid formulas [24, Theorem 4]; for concreteness' sake we set $\mathcal{G} = (I, \{1\})$.

Consider some axiomatization of the intuitionistic propositional calculus IPC, such as Kleene's [52]:

Ax1. $\alpha \to (\beta \to \alpha)$;

Ax2. $(\alpha \to \beta) \to ((\alpha \to (\beta \to \gamma)) \to (\alpha \to \gamma))$;

Ax3. $\alpha \to (\beta \to \alpha \wedge \beta)$;

Ax4. $\alpha \wedge \beta \to \alpha$;

Ax5. $\alpha \wedge \beta \to \beta$;

Ax6. $\alpha \to \alpha \vee \beta$;

Ax7. $\beta \to \alpha \vee \beta$;

Ax8. $(\alpha \to \gamma) \to ((\beta \to \gamma) \to (\alpha \vee \beta \to \gamma))$;

Ax9. $(\alpha \to \beta) \to ((\alpha \to \neg\beta) \to \neg\alpha)$;

Ax10. $\neg\alpha \to (\alpha \to \beta)$.

Any theorem of IPC is valid in all the $\mathcal{G}_n$'s, as well as in $\mathcal{G}$. On the other hand, Gödel proved in [34] that no finite matrix is weakly adequate for the consequence operation $C_{\mathrm{ICP}}$ ($C_{\mathrm{ICP}}$ is defined thus: $A \in C_{\mathrm{ICP}}(\Gamma)$ iff $A$ is deducible from $\Gamma$ via Ax1–Ax10 and modus ponens). Gödel's proof can also be found as an exercise in [17, p. 145]. In [24], Dummett characterizes the set of formulas valid in $\mathcal{G}$ (equivalently, valid in every $\mathcal{G}_n$) as the set of formulas deducible from Ax1–Ax11, where Ax11 is the axiom schema expressing linearity:

Ax11. $(\alpha \to \beta) \vee (\beta \to \alpha)$.

In [48], Jáskowski constructs an infinite family of matrices whose corresponding sets of valid formulas converge to $C_{\mathrm{ICP}}(\emptyset)$. This is a classical example of an interesting phenomenon: many finitely axiomatizable propositional logics can be approximated by sequences of finite-valued logics (see [4]). Jáskowski constructs also an infinite matrix (different from the Lindenbaum algebra of IPC) whose set of valid formulas coincides with the set of intuitionistic theorems.

We remark that by a result of Wroński's [112], no denumerable matrix is strongly adequate to $C_{\mathrm{ICP}}$.

### 2.5.3 Belnap's system

Belnap's logic is a four-valued non-linearly-ordered logic [6], [49]. It is intended to deal with incomplete and possibly contradictory knowledge bases. Suppose a computer is designed to argue about a collection of ground facts. Suppose that it is informed about the truth or falseness of these facts by different agents, each of them globally reliable, but sometimes incorrect in some specific statement. It may happen that the computer is told first that the fact $p$ is true, and then that it is false; this kind of exposure to contradiction is common in our everyday's life. A perfect Boolean reasoner that faces a contradiction breaks down, and starts inferring as a new theorem every statement he can think of; of course, this is not what we do, nor what we want our computer do.

Let **4** be the powerset of **2**; for simplicity's sake, we set $\mathbf{4} = \{\emptyset, 0, 1, 01\}$. If $p$ is a ground fact, i.e., a propositional variable, then four possibilities arise:

- either the computer has been told that $p$ is true;
- or that it is false;
- or that it is true *and* false (say by two different agents, or in different circumstances);
- or the computer has not been told anything.

These four possibilities correspond to assignments to $p$ of values from **4** (the values 1, 0, 01, and $\emptyset$, respectively). For each $p$, its successive truth-assignments must be increasing with respect to the order given by the *approximation lattice*

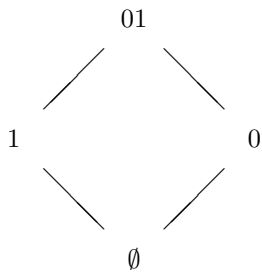

On the other hand, reasoning about compound facts is made with respect to the *logical lattice*

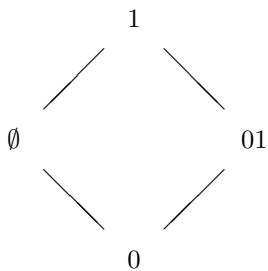

so that the tables for negation, disjunction, and conjunction are as follows:

| $\neg$ | |
|---|---|
| $\emptyset$ | $\emptyset$ |
| 0 | 1 |
| 1 | 0 |
| 01 | 01 |

| $\vee$ | $\emptyset$ | 0 | 1 | 01 |
|---|---|---|---|---|
| $\emptyset$ | $\emptyset$ | $\emptyset$ | 1 | 1 |
| 0 | $\emptyset$ | 0 | 1 | 01 |
| 1 | 1 | 1 | 1 | 1 |
| 01 | 1 | 01 | 1 | 01 |

| $\wedge$ | $\emptyset$ | 0 | 1 | 01 |
|---|---|---|---|---|
| $\emptyset$ | $\emptyset$ | 0 | $\emptyset$ | 0 |
| 0 | 0 | 0 | 0 | 0 |
| 1 | $\emptyset$ | 0 | 1 | 01 |
| 01 | 0 | 0 | 01 | 01 |

The key feature of Belnap's logic lies in the interplay between the approximation and the logical lattices [28], [29], [33]. There are no distinguished values, nor an implication connective. Entailment is metalinguistic: we say that $A$ entails $B$ if, for every structure $\sigma$, we have $\sigma(A) \leq \sigma(B)$ in the logical lattice. In this case, we write $A \leq B$; we write $A \equiv B$ if both $A \leq B$ and $B \leq A$ hold (i.e., $A$ and $B$ are semantically equivalent). As $\vee$ and $\wedge$ are defined in the standard way with respect to the logical lattice, it follows that all Boolean equivalences between formulas not involving $\neg$ still hold in Belnap's logic. In particular, $\vee$ and $\wedge$ are associative, commutative, and mutually distributive. One easily sees that $\neg\neg A \equiv A$ and that the De Morgan's laws hold. On the other hand, neither *tertium non datur* $A \vee \neg A \equiv 1$ nor *ex falso quodlibet* $A \wedge \neg A \leq B$ are generally valid; of course, this is welcome in a system which is to deal with contradiction.

## 2.6 Extension to the predicate case

It is only fair to say that, if a unified framework for propositional many-valued logic is still to come, much more so is the situation for the predicate version. Beyond the serious technical problems, what is missing is a convincing philosophical and linguistic interpretation. In most cases, authors simply skip any discussion about "meaning", and go directly to the technical content.

**Definition 2.6.1** A language $L$ for predicate many-valued logic is given by:

1. a —non-empty— finite or denumerable set of predicate symbols $P, R, U, \ldots$, each one with a fixed arity $\geq 0$;

2. a —possibly empty— finite or denumerable set of function symbols $f, g, h, \ldots$, each one with a fixed arity $\geq 0$;

3. a finite set of connectives $c_1, \ldots, c_m$, each one with a fixed arity $\geq 0$;

4. a finite set of quantifiers $Q_1, \ldots, Q_s$;

5. a denumerable set of individual variables $x_1, x_2, x_3, \ldots$

A *matrix* for $L$ is a system $\mathcal{M} = (M, D, q_1, \ldots, q_s)$, where $(M, D)$ is a matrix for the propositional fragment of $L$, and each $q_i$ is a function : $\mathcal{P}(M) \setminus \{\emptyset\} \to M$. By the propositional fragment of $L$, we mean the language $L$ with (2), (4),
(5) above dropped, and all predicate symbols of arity 0.
A *structure* for $(L, \mathcal{M})$ is a pair $\mathcal{S} = (S, \sigma)$, where:

1. $S$ is a non-empty set;

2. $\sigma$ is a mapping that associates:

   - to each function symbol $f$ of arity $u$, a function $\sigma(f) : S^u \to S$;
   - to each predicate symbol $P$ of arity $u$, a mapping $\sigma(P) : S^u \to M$.

In other words, any function symbol is realized by $\sigma$ as an actual function over $S$, while every predicate symbol is realized as a fuzzy subset of $S$. 0ary predicate symbols are propositional symbols, and 0ary function symbols are realized as individual constants (i.e., elements of $S$). We remark that the language does not contain a distinguished predicate symbol = for equality of individuals.

Let us denote by $L(S)$ the language $L$ enriched with individual constants, one for each element of $S$, and let $SENT(L), SENT(L(S))$ be the subsets of $FORM(L), FORM(L(S))$, respectively, consisting of the formulas with no free individual variables. Then $\mathcal{S}$ induces a mapping $\mathcal{S} : SENT(L(S)) \to M$ as follows:

- if $P$ is a predicate symbol of arity $u$ and $t_1, \ldots, t_u$ are closed terms, then $\mathcal{S}(P(t_1, \ldots, t_u)) = (\sigma(P))(\sigma(t_1), \ldots, \sigma(t_u))$;

- if $c$ is a connective of arity $u$ and $A_1, \ldots, A_u \in SENT(L(S))$, then $\mathcal{S}(cA_1 \ldots A_u) = c(\mathcal{S}(A_1), \ldots, \mathcal{S}(A_u))$;

- if $Q$ is a quantifier whose corresponding function : $\mathcal{P}(M) \setminus \{\emptyset\} \to M$ is $q$, then $\mathcal{S}(QxA(x)) = q(\{\mathcal{S}(A(a)) : a \in S\})$.

The set $\{\mathcal{S}(A(a)) : a \in S\}$ is the *distribution* of $A(x)$ in $\mathcal{S}$, hence the name *distribution quantifiers* [64], [13], [14].

A formula $A(x_1, \ldots, x_n) \in FORM(L)$ containing the free individual variables $x_1, \ldots, x_n$, is *satisfied* in $\mathcal{S}$ if there exist $a_1, \ldots, a_n \in S$ such that $\mathcal{S}(A(a_1, \ldots, a_n)) \in D$, and it is *true* in $\mathcal{S}$ if, for every $a_1, \ldots, a_u \in S$, $\mathcal{S}(A(a_1, \ldots, a_u)) \in D$. In classical logic, where the universal quantifier $\forall$ is available, $A$ is true in $\mathcal{S}$ iff so is its universal closure; in our case, one must resort to the above definition.

In a given propositional matrix $\mathcal{M} = (M, D)$ of cardinality $n$, one can define $n^{2^n-1}$ different quantifier functions. An interesting and difficult problem is to determine how many of them are needed to define all the others. For example, in the two-element Boolean algebra $\mathbf{2}$, the function

$$q_\forall : \begin{cases} \{0\} \mapsto 0; \\ \{1\} \mapsto 1; \\ \{0,1\} \mapsto 0; \end{cases}$$

together with the connectives $\vee, \neg$, defines all the eight possible quantifier functions.

In classical logic, the universal quantifier is naturally associated with the connective $\wedge$. Explicitly, $\wedge$ induces an inf-semilattice structure on $\mathbf{2}$. As $\mathbf{2}$ is finite, it results to be a complete inf-semilattice, and the $q_\forall$-image of a nonempty subset of $\mathbf{2}$ is the g.l.b. of the subset; analogous considerations hold for the existential quantifier $\exists$ and the connective $\vee$. In general, the situation is as follows: any binary, idempotent, associative and commutative connective $c$ defined in the (possibly infinite) algebra $M$ induces a partial order by $a \leq b$ iff $c(a, b) = a$; this order makes $M$ an inf-semilattice. Conversely, each semilattice structure over $M$ determines a binary, idempotent, associative and commutative connective. If our semilattice happens to be complete (as it is always the case if $M$ is finite), then it yields a quantifier function $q_c : \mathcal{P}(M) \setminus \{\emptyset\} \to M$, defined by $q_c(X) = \inf X$. If $M$ has cardinality $n$, then there are $n^{n-1}$ possible semilattice structures over $M$ [114, §1.7], and so $n^{n-1}$ possible quantifiers of the above form.

In order to extend sequent and tableau calculi, we need introduction —respectively, elimination— rules for quantifiers. We refer to [92] for the extension of resolution techniques to the predicate case. Let $\mathcal{S} = (S, \sigma)$ be a structure. Define the class $\Phi(S)$ of *quantified signed formula expressions with parameters in $S$* as follows:

- for any $A \in FORM(L(S))$ and any $i \in M$, $A^i \in \Phi(S)$;

- if $E, F \in \Phi(S)$, then $E \wedge F, E \vee F \in \Phi(S)$;

- if $E \in \Phi(S)$ and $x$ is an individual variable, then $\forall xE, \exists xE \in \Phi(S)$.

$\mathcal{S}$ induces a mapping $\overline{\mathcal{S}}$ from the set of elements of $\Phi(S)$ having no free variables to the two-element Boolean algebra by induction:

- if $A \in SENT(L(S))$, then $\overline{\mathcal{S}}(A^i) = 1$ iff $\mathcal{S}(A) = i$;

- $\overline{\mathcal{S}}$ distributes over $\wedge$ and $\vee$;

- $\overline{\mathcal{S}}(\forall xEx) = \min\{\overline{\mathcal{S}}(Ea) : a \in S\}$ and $\overline{\mathcal{S}}(\exists xEx) = \max\{\overline{\mathcal{S}}(Ea) : a \in S\}$.

We denote by $\Phi$ the class of all quantified signed formula expressions having no free variables and no parameters from any $S$. An element of $\Phi$ is in CNF (DNF) if it is of the form $\forall \bar{u} \exists \bar{v} E(\bar{u}, \bar{v})$ ($\exists \bar{v} \forall \bar{u} E(\bar{u}, \bar{v})$), where $\bar{u}, \bar{v}$ are strings of variables and $E(\bar{u}, \bar{v})$ is quantifier-free and in conjunctive (disjunctive) normal form.

If $E, F \in \Phi$ are such that, for any structure $\mathcal{S}$, $\overline{\mathcal{S}}(E) = \overline{\mathcal{S}}(F)$, then we say that $E, F$ are *equivalent*, and we write $E \equiv F$. As in the propositional case, if $E, F$ are equivalent in the classical predicate calculus, then they are equivalent as quantified signed formula expressions.

Let us fix a unary predicate symbol $P$. Let $Q$ be a quantifier whose associated function is $q$, and let $i \in M$. We define an *$i$th CNF (DNF) for $Q$* to be some $E \in \Phi$ such that:

- $E$ is in CNF (DNF);
- the only signed formulas appearing in $E$ are of the form $(Pw)^j$, for some variable $w$ and truth-value $j$ (of course, $w$ gets quantified in the prefix of $E$);
- $E \equiv (QxPx)^i$.

Usually it is not difficult to construct an $i$th CNF and an $i$th DNF for a quantifier, given the corresponding quantifier function $q$. We refer to [114], [40, §3.3, §5.4] for algorithms for systematically constructing such normal forms.

Suppose $\forall u_1 \ldots \forall u_p \exists v_1 \ldots \exists v_q (D_1 \wedge \cdots \wedge D_h)$ is an $i$th CNF for $Q$. Then we have the following introduction rule for $Q$ at place $i$:

$$\frac{G, F_1 \quad G, F_2 \quad \cdots \quad G, F_h}{G, (QxAx)^i}$$

where $G$ is an arbitrary sequent, and $F_s$ is the sequent obtained from $D_s$ by replacing $(Pu_e)^j$ by $(Az_e)^j$ and $(Pv_f)^j$ by $(At_f)^j$. Here $z_1, \ldots, z_p$ are variables which do not appear anywhere in $G, (QxAx)^i$ (the *eigenvariables* of the rule), while $t_1, \ldots, t_q$ are arbitrary terms.

**Example 2.6.2** Consider the universal quantifiers $\forall, \exists$ in classical logic. Then
$(\forall x Px)^0 \equiv \exists v (Pv)^0$ and $(\exists x Px)^0 \equiv \forall u (Pu)^0$. Hence the familiar rules

$$\frac{\Gamma, At \Rightarrow \Delta}{\Gamma, \forall x Ax \Rightarrow \Delta} \qquad \frac{\Gamma, Az \Rightarrow \Delta}{\Gamma, \exists x Ax \Rightarrow \Delta}$$

where $t$ is an arbitrary term and $z$ is a variable that does not appear in $\Gamma, \exists x Ax \Rightarrow \Delta$.

On the other hand, let $\exists v_1 \ldots \exists v_q \forall u_1 \ldots \forall u_p (C_1 \vee \cdots \vee C_h)$ be an $i$th DNF for $Q$. We have the tableau rule

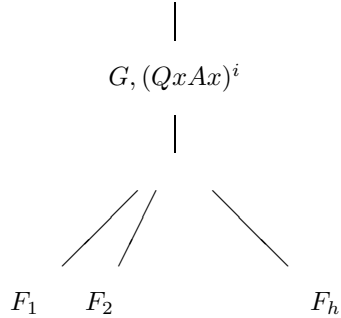

where $F_s$ is obtained from $C_s$ by substituting $(Pv_f)^j$ with $(Az_f)^j$ and $(Pu_e)^j$ with $(At_e)^j$. Here $z_1, \ldots, z_q$ are fresh variables, nowhere appearing in the branches through $G, (QxAx)^i$, and $t_1, \ldots, t_p$ are arbitrary terms. $(QxAx)^i$ is declared to be used only if our DNF does not contain universal quantifiers (i.e., there are no terms $t_1, \ldots, t_p$). Alternatively, if one looks for a systematic way to construct the tableau, one may declare $(QxAx)^i$ used, and add $(QxAx)^i$ to the end of the new branches [89], [13].

**Example 2.6.3** We introduce the quantifiers $\forall$ and $\exists$ in Łukasiewicz logic as generalized $\wedge$ and $\vee$, respectively. For $\emptyset \neq X \subseteq I$, we define $q_\forall(X) = \inf X$ and $q_\exists(X) = \sup X$. Let us consider the three-valued version. We have $(\forall x Px)^{1/2} \equiv \exists v (Pv)^{1/2} \wedge \forall u ((Pu)^{1/2} \vee (Pu)^1)$. It is easy to see that in the classical predicate calculus this last formula is equivalent to both $\exists v \forall u (((Pv)^{1/2} \wedge (Pu)^{1/2}) \vee ((Pv)^{1/2} \wedge (Pu)^1))$ and $\exists v \forall u ((Pu)^{1/2} \vee ((Pv)^{1/2} \wedge (Pu)^1))$. Hence we can choose between the rules

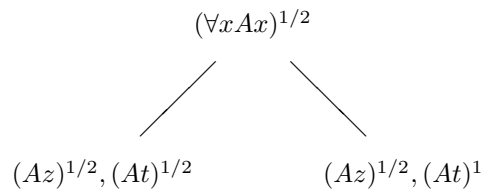

and

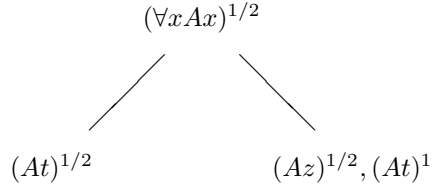

In both cases, $z$ is a fresh variable, $t$ is an arbitrary term, and the occurrence of $(\forall xAx)^{1/2}$ is not declared to be used.

In the infinite-valued Łukasiewicz logic, we do not have analytic systems at our disposal, and we must resort to Hilbert systems. This does not carry us much further. Indeed, it is a classical result that the set of predicate formulas valid in $I$ is not recursively enumerable [85], [79]. On the other hand, the set of formulas whose value is greater than 0 in any structure has been axiomatized in [47].

## 3 MV-algebras

As we said in the Introduction, this part is more mathematically oriented, and presupposes some acquaintance with universal algebra.

Beyond its intrinsic interest, the algebrization of a logic yields useful techniques. First of all, many completeness proofs are based on representation theorems for certain free algebras. This is the case for the three logics of Example 2.4.2, and in this section we sketch the proof for Łukasiewicz logic. Secondly, one can look for normal forms for formulas. These forms are admittedly scarce, but we already encountered one such example at the end of §2.3, and in the case of Łukasiewicz logic we cite the results obtained in [69], [75]. Thirdly, one is allowed to see proofs as computations in some algebra: this is apparent, e.g., in Stachniak's approach to resolution. [80] and [10] are basic references for algebraic logic.

In this section we discuss the algebrization of Łukasiewicz logic, and study the corresponding algebras (MV-algebras). MV-algebras stand to Łukasiewicz logic as Boolean algebras stand to classical logic. They form an interesting class of mathematical objects, bearing relevant connections with areas of classical mathematics such as ordered group theory, functional analysis, and algebraic geometry. The book [19] is a systematic account of basic results otherwise scattered in the literature.

### 3.1 Basic facts

We recall that in Łukasiewicz logic the sets of connectives $(\to, \neg, 1)$, $(\to, 0)$, $(\oplus, \otimes, \neg, 0, 1)$ are interdefinable. For all this section we fix $L$ to be the language $(\oplus, \otimes, \neg, 0, 1)$, of similarity type $2, 2, 1, 0, 0$. We denote by $I$ the abstract algebra obtained from the matrix $\mathcal{L} = (I, \{1\})$ of §2.4 by dropping reference to the set $\{1\}$ of designated values, and by changing the language from $(\to, \neg, 1)$ to $L$. To be more specific, $I$ has as base set the set of all real numbers between 0 and 1, and the operations are defined thus:

$$
\begin{aligned}
i \oplus j &= \min\{i+j, 1\} \\
i \otimes j &= \max\{i+j-1, 0\} \\
\neg i &= 1 - i \\
0 &= \text{the number } 0 \\
1 &= \text{the number } 1
\end{aligned}
$$

Let us rewrite the Łukasiewicz axioms of §2.4 as identities:

Ł1. $\neg a \oplus (\neg b \oplus a) = 1$;

Ł2. $\neg(\neg a \oplus b) \oplus (\neg(\neg b \oplus c) \oplus (\neg a \oplus c)) = 1$;

Ł3. $\neg(\neg\neg a \oplus \neg b) \oplus (\neg b \oplus a) = 1$;

Ł4. $\neg(\neg(\neg a \oplus b) \oplus b) \oplus (\neg(\neg b \oplus a) \oplus a) = 1$.

A bit of computation shows that Ł1–Ł4 are valid in $I$.

**Definition 3.1.1** An MV-*algebra* is an abstract algebra $N = (N, \oplus, \otimes, \neg, 0, 1)$ satisfying the identities Ł1–Ł4. We denote by **MV** the equational class of all MV-algebras.

MV-algebras have been introduced by Chang in [15]. MV stands for *Many Valued*; this name is well established in the literature, even though these algebras refer to one specific many valued logic, namely the Lukasiewicz one. Setting $a \vee b = (a \to b) \to b = \neg(\neg a \oplus b) \oplus b$ and $a \wedge b = \neg(\neg a \otimes b) \otimes b$, we have that $(N, \vee, \wedge, 0, 1)$ is a bounded distributive lattice [15], [16]. As usual, we write $a \leq b$ if $a \vee b = b$.

Axioms Ł1–Ł4 are not convenient to work with. Chang's original axioms are the following:

C1. $a \oplus b = b \oplus a$; $\qquad\qquad\qquad a \otimes b = b \otimes a$;

C2. $a \oplus (b \oplus c) = (a \oplus b) \oplus c$; $\qquad a \otimes (b \otimes c) = (a \otimes b) \otimes c$;

C3. $a \oplus \neg a = 1$; $\qquad\qquad\qquad a \otimes \neg a = 0$;

C4. $a \oplus 1 = 1$; $\qquad\qquad\qquad a \otimes 0 = 0$;

C5. $a \oplus 0 = a$; $\qquad\qquad\qquad a \otimes 1 = a$;

C6. $\neg(a \oplus b) = \neg a \otimes \neg b$; $\qquad\qquad \neg(a \otimes b) = \neg a \oplus \neg b$;

C7. $\neg\neg a = a$;

C8. $\neg 0 = 1$;

C9. $a \vee b = b \vee a$; $\qquad\qquad\qquad a \wedge b = b \wedge a$;

C10. $a \vee (b \vee c) = (a \vee b) \vee c$; $\qquad a \wedge (b \wedge c) = (a \wedge b) \wedge c$;

C11. $a \oplus (b \wedge c) = (a \oplus b) \wedge (a \oplus c)$; $\qquad a \otimes (b \vee c) = (a \otimes b) \vee (a \otimes c)$.

A more compact axiomatization is due to Mangani and Mundici [65]:

M1. $a \oplus b = b \oplus a$;

M2. $a \oplus (b \oplus c) = (a \oplus b) \oplus c$;

M3. $a \oplus 0 = a$;

M4. $a \oplus 1 = 1$;

M5. $\neg\neg a = a$;

M6. $\neg 0 = 1$;

M7. $\neg(\neg a \oplus b) \oplus b = \neg(\neg b \oplus a) \oplus a$;

M8. $a \otimes b = \neg(\neg a \oplus \neg b)$.

In every MV-algebra, the following hold:

1. De Morgan's laws for $\vee, \wedge, \neg$;

2. $a \otimes b \leq a \wedge b \leq a \leq a \vee b \leq a \oplus b$;

3. $(a \to b) \vee (b \to a) = 1$;

4. $a \vee b = b$ iff $\neg a \oplus b = 1$;

5. $a \oplus b = b$ iff $\neg a \vee b = 1$;

6. $a \oplus (b \vee c) = (a \oplus b) \vee (a \oplus c)$;

7. $a \otimes (b \wedge c) = (a \otimes b) \wedge (a \otimes c)$.

As $\oplus, \otimes, \vee, \wedge$ are associative and commutative, we freely use expressions of the form $a_1 \oplus a_2 \oplus \cdots \oplus a_k$ and the like. We set $0a = 0$, $(k+1)a = a \oplus ka$ and $a^0 = 1$, $a^{k+1} = a \otimes a^k$. We say that $a$ has *order* $k$ if $k$ is the least positive integer with $ka = 1$; if no such $k$ exists, we say that $a$ has infinite order.

**Example 3.1.2** We already introduced the MV-algebra $I$. Let us denote by $[0,1]^n$ the set of all $n$-uples of elements of $[0,1]$; $[0,1]^\omega$ is the set of all countable sequences of elements of $[0,1]$. A *McNaughton function over* $[0,1]^n$ is a continuous function $f : [0,1]^n \to \mathbf{R}$ for which the following condition holds:

there exist linear affine polynomials $f_1, \ldots, f_k$, each $f_i$ of the form $a_i^1 x_1 + a_i^2 x_2 + \cdots + a_i^n x_n + a_i^{n+1}$ for certain integers $a_i^1, \ldots, a_i^{n+1}$, such that, for every $\sigma \in [0,1]^n$, there exists $1 \leq i \leq k$ with $f(\sigma) = f_i(\sigma)$.

A *McNaughton function over* $[0,1]^\omega$ is a function $g : [0,1]^\omega \to \mathbf{R}$ such that there exists a McNaughton function $f$ over $[0,1]^n$ for which $g = f \circ \rho_n$ ($\rho_n$ is the natural projection : $[0,1]^\omega \to [0,1]^n$). We let $M_n$ be the subset of $I^{[0,1]^n}$ whose elements are McNaughton functions over $[0,1]^n$; it is easy to see that $M_n$ is in fact a subalgebra, hence an MV-algebra. Analogously, $M_\omega$ is the MV-algebra of all McNaughton functions over $[0,1]^\omega$ whose range is in $[0,1]$. We let $\mathbf{0}$ and $\mathbf{1}$ be the functions, of domain $[0,1]^n$ or $[0,1]^\omega$ according to the context, whose values are identically 0 or 1, respectively.

Let $FORM_n(L)$ be the subset of $FORM(L)$ whose elements are built up from the propositional variables $p_1, \ldots, p_n$ only. We identify $FORM(L)$ with the set of terms in the language of MV-algebras; when we are thinking of $A \in FORM_n(L)$ as a term, we substitute the free variables $x_1, \ldots, x_n$ for the propositional variables $p_1, \ldots, p_n$. We also identify the space of all structures with $[0,1]^\omega$. If $A$ is a formula, then $A$ induces a mapping $\mathbf{A} : [0,1]^\omega \to [0,1]$ in the natural way: $\mathbf{A}(\sigma) = \sigma(A)$. If $A \in FORM_n(L)$, then the value of $\sigma(A)$ is unaffected by changing the value of $\sigma(q)$, for $q \notin \{p_1, \ldots, p_n\}$. It follows that, with respect to $A$, the space of all structures is identifiable with the $n$-cube, and we may regard $\mathbf{A}$ as a function $\mathbf{A} : [0,1]^n \to [0,1]$. An easy induction over the complexity of $A$ shows that $\mathbf{A}$ is a McNaughton function, hence an element of $M_\omega$ or $M_n$, respectively.

For $A, B \in FORM(L)$, we write $A \equiv B$ if the identity of terms $A = B$ holds in every MV-algebra; equivalently, the identity $A = B$ may be derived from Ł1–Ł4,C1–C11,M1–M8 by using the rules of the equational calculus. The relation $\equiv$ is a congruence, and the algebras $F_\omega = FORM(L)/\equiv$ and $F_n = FORM(L)_n/\equiv$ are the free MV-algebras over $\omega$ and $n$ generators, respectively. The function $(A/\equiv) \mapsto \mathbf{A}$ is a homomorphism of MV-algebras, either from $F_\omega$ to $M_\omega$, or from $F_n$ to $M_n$.

Let us write $C_{\text{LPC}}$ for the consequence operation induced by the axiom schemas Ax1–Ax4 of §2.4, along with modus ponens. The completeness theorem for the Łukasiewicz infinite-valued calculus can be expressed in any of the following forms:

1. $\mathcal{L}$ is weakly adequate for $C_{\text{LPC}}$;

2. $A \equiv 1$ iff $\mathbf{A} = \mathbf{1}$;

3. for every $n$, the homomorphism $(A/\equiv) \mapsto \mathbf{A} : F_n \to M_n$ is injective;

4. the homomorphism $(A/\equiv) \mapsto \mathbf{A} : F_\omega \to M_\omega$ is injective;

5. $A \equiv B$ iff $\mathbf{A} = \mathbf{B}$;

6. $F_\omega$ belongs to the subvariety of $\mathbf{MV}$ generated by $I$;

7. $I$ generates $\mathbf{MV}$.

As remarked in §2.4, there are several different proofs of these statements. We sketch Chang's proof of (6) [16]. As a first step, Chang shows that every MV-algebra is a subdirect product of linearly ordered MV-algebras (we will see this in §3.2). This reduces the problem to proving that every totally ordered MV-algebra belongs to the subvariety of $\mathbf{MV}$ generated by $I$, i.e., that every identity which is true in $I$ is also true in every totally ordered MV-algebra.

**Definition 3.1.3** A *lattice ordered abelian group* ($\ell$-group) is an algebra $(G, +, -, 0, \vee, \wedge)$ such that $(G, +, -, 0)$ is an abelian group, $(G, \vee, \wedge)$ is a lattice, and the identity $a + (b \vee c) = (a+b) \vee (a+c)$ holds. A *totally ordered abelian group* (*o*-group) is an $\ell$-group in which the lattice order is total. A *strong unit* of the $\ell$-group $G$ is an element $u > 0$ of $G$ such that, for every $a \in G$, there exists $m \in \mathbf{N}$ with $a \leq mu$.

It is well known [107] that:

a. every universal formula in the language $(+, -, 0, \vee, \wedge)$ which is true in $\mathbf{Q}$ is true in all *o*-groups.

Chang's fundamental idea is to associate:

b. to every totally ordered MV-algebra $N$ a pair $(G_N, u_N)$, where $G_N$ is an o-group, and $u_N$ is a strong unit of $G_N$; in particular, $(\mathbf{R}, 1)$ corresponds to $I$;

c. to every identity $A = B$ in the language of MV-algebras with free variables $x_1, \ldots, x_n$, a quantifier-free formula $\varphi_{AB}(x_1, \ldots, x_n, y)$ in the language of $\ell$-groups;

in such a way that:

d. $N \models \forall x_1 \ldots \forall x_n (A = B)$ iff $G_N \models \forall x_1 \ldots \forall x_n \varphi_{AB}(x_1, \ldots, x_n, u_N)$.

Let now $A = B$ be any identity in the language of MV-algebras with free variables $x_1, \ldots, x_n$, let $N$ be a totally ordered MV-algebra, and assume $I \models \forall x_1 \ldots \forall x_n (A = B)$. We then have:

$\mathbf{R} \models \forall x_1 \ldots \forall x_n \varphi_{AB}(x_1, \ldots, x_n, 1)$, by (b) e (d);

$\mathbf{Q} \models \forall x_1 \ldots \forall x_n \varphi_{AB}(x_1, \ldots, x_n, 1)$, because $\mathbf{Q}$ is an $\ell$-subgroup of $\mathbf{R}$;

$\mathbf{Q} \models \forall x_1 \ldots \forall x_n \forall y (0 < y \to \varphi_{AB}(x_1, \ldots, x_n, y))$, because, for every $0 < c \in \mathbf{Q}$, there exists an automorphism of $\mathbf{Q}$ that maps 1 in $c$;

$G_N \models \forall x_1 \ldots \forall x_n \forall y (0 < y \to \varphi_{AB}(x_1, \ldots, x_n, y))$, by (a);

$G_N \models \forall x_1 \ldots \forall x_n \varphi_{AB}(x_1, \ldots, x_n, u_N)$;

$N \models \forall x_1 \ldots \forall x_n (A = B)$, by (d).

It follows that the validity of $A = B$ in $I$ implies the validity of $A = B$ in $N$, which is what we wanted to prove.

In force of the completeness theorem, we can test if an identity holds in every MV-algebra by testing it in $I$. This fits nicely with the situation for Boolean algebras, where we test identities in $\mathbf{2}$. As McNaughton functions are continuous, we can actually test identities in the subalgebra of $I$ whose elements are the rational numbers between 0 and 1. It is a fortunate case that the computational complexity of the tautology problem does not increase: the set $C_{\text{LPC}}(\emptyset) = C_{\mathcal{L}}(\emptyset)$ is co-NP complete [66].

Again by completeness, we know that $F_n$ and $F_\omega$ are subalgebras of $M_n$ and $M_\omega$. Actually, something stronger happens: McNaughton's theorem [59] is to the effect that the embedding $(A/\equiv) \mapsto \mathbf{A}$ is surjective, so that $M_n$ and $M_\omega$ are isomorphic to the free MV-algebras over $n$ and $\omega$ generators. This allows us to freely identify the formula $A$ with its equivalence class $(A/\equiv)$, and with the function $\mathbf{A}$; we shall do this in the sequel.

In [65], Mundici extends the correspondence $N \mapsto (G_N, u_N)$ to arbitrary MV-algebras.

**Definition 3.1.4** Let $G, H$ be $\ell$-groups, $u$ a fixed strong unit of $G$, and $v$ a fixed strong unit of $H$. A *unital homomorphism* $\varphi : (G, u) \to (H, v)$ is an $\ell$-group homomorphism $\varphi : G \to H$ such that $\varphi(u) = v$. The *category of $\ell$-groups with strong unit* is the category whose objects are all pairs $(G, u)$, where $G$ is an $\ell$-group and $u$ is a strong unit of $G$, and whose morphisms are the unital homomorphisms.

**Proposition 3.1.5** [65, Theorem 2.5] *Let $(G, u)$ be an $\ell$-group with strong unit $u$. Let $\Gamma(G, u)$ be the structure*
$$\Gamma(G, u) = ([0, u], \oplus, \otimes, \neg, 0, 1)$$
*defined as follows:*

$$\begin{aligned}
[0, u] &= \{a \in G : 0 \leq a \leq u\} \\
a \oplus b &= (a + b) \wedge u \\
a \otimes b &= (a + b - u) \vee 0 \\
\neg a &= u - a \\
0 &= \text{the additive identity 0 of } G \\
1 &= u
\end{aligned}$$

*Then $\Gamma(G, u)$ is an MV-algebra, and the lattice order induced by the MV operations coincides with the order inherited from the $\ell$-group $G$. Let $\varphi : (G, u) \to (H, v)$ be a unital homomorphism, and let $\Gamma \varphi$ be the restriction of $\varphi$ to $[0, u]$. Then the image of $\Gamma \varphi$ is contained in $[0, v]$, and $\Gamma \varphi : \Gamma(G, u) \to \Gamma(H, v)$ is a homomorphism of MV-algebras.* ∎

**Example 3.1.6** The set $M[0,1]^n$ of all McNaughton functions over $[0,1]^n$ is an $\ell$-subgroup of the $\ell$-group $\mathbf{R}^{[0,1]^n}$, and $\mathbf{1}$ is a strong unit of it. We have $M_n = \Gamma(M[0,1]^n, \mathbf{1})$. Analogously, the set $M[0,1]^\omega$ of all McNaughton functions over $[0,1]^\omega$ is an $\ell$-group and $M_\omega = \Gamma(M[0,1]^\omega, \mathbf{1})$.

**Example 3.1.7** Let $\mathbf{Z} \oplus_{\text{lex}} \mathbf{Z}$ be the $o$-group whose underlying abelian group is the direct sum of two copies of $\mathbf{Z}$, and whose order is lexicographic:

$$(a,b) \leq (c,d) \text{ iff } (a < c \text{ or } (a = c \text{ and } b \leq d)).$$

*Chang's algebra* is the MV-algebra $C = \Gamma(\mathbf{Z} \oplus_{\text{lex}} \mathbf{Z}, (1,0))$ [15, p. 474]. The order type of $C$ consists of two copies of the natural numbers, each facing the other:

$$(0,0) < (0,1) < (0,2) < \cdots < (1,-2) < (1,-1) < (1,0).$$

$C$ is the prototypical example of a linearly ordered nonsimple MV-algebra (see §3.2). The element $(0,1)$ generates $C$, which is therefore a quotient of $M_1$.

**Theorem 3.1.8** [65, Theorem 3.9] $\Gamma$ *is a full, faithful, and representative functor (i.e., a categorical equivalence) between the category of $\ell$-groups with strong unit and the category of MV-algebras. In particular, for every MV-algebra $N$, there exists a unique $\ell$-group with strong unit $(G_N, u_N)$ such that $N$ is isomorphic to $\Gamma(G_N, u_N)$. If $N$ is countable, then $G_N$ is countable.* ∎

The $\Gamma$ functor allows us to translate facts and problems from MV-algebras to $\ell$-groups, and conversely. The techniques in one field may be exported to the other and, since $\ell$-groups have relevant connections to analysis (indeed, typical $\ell$-groups are groups of continuous functions) this correspondence can be further extended. As an example of this interplay, we cite [70], where all $C^*$-algebras whose $K_0$ group is lattice ordered, and which are limits of sequences of finite-dimensional $C^*$-algebras, have been classified by means of MV-algebras.

## 3.2 Ideals and representations

Let $N$ be an MV-algebra. An *ideal* of $N$ is a subset $J$ of $N$ containing 0, closed under $\oplus$, and such that $a \leq b \in J$ implies $a \in J$; $J$ is *proper* if $J \neq M$. Dually, a *filter* of $N$ is a subset $F$ of $N$ containing 1, closed under $\otimes$, and such that $a \geq b \in F$ implies $a \in F$. Of course, $J$ is an ideal iff $\{\neg a : a \in J\}$ is a filter, and conversely. It is also easy to see that $F$ is a filter iff $F$ is closed under modus ponens ($a, a \to b \in F$ implies $b \in F$) and contains 1; hence filters correspond to theories in the infinite-valued calculus. As usual, every ideal $J$ induces the congruence $\sim_J$ by

$$a \sim_J b \quad \text{iff} \quad \neg(a \leftrightarrow b) = (a \otimes \neg b) \vee (b \otimes \neg a) \in J$$

and the equivalence class of 0 under any congruence is an ideal. We write $N/J$ for the image of $N$ under the natural epimorphism $a \mapsto a/J = a/\sim_J$.

If $X$ is any subset of $N$, then the filter $\overline{X}$ generated by $X$ is the intersection of all filters containing $X$. It is not difficult to see that $\overline{X}$ coincides with the set of all $a \in N$ for which there exist $b_1, b_2, \ldots, b_r \in X$ (not necessarily distinct) with $b_1 \otimes b_2 \otimes \cdots \otimes b_r \leq a$. If we translate this in terms of the Łukasiewicz calculus, we obtain the infinite-valued version of the Deduction Theorem [82, Theorem 2.2]:

$$A_1, \ldots, A_r \vdash B \quad \text{iff} \quad \vdash (A_1 \otimes A_2 \otimes \cdots \otimes A_r)^t \to B \text{ for some } t \geq 1.$$

Let $J$ be a proper ideal of the MV-algebra $N$. $J$ is *maximal* if it is not properly contained in any proper ideal. $J$ is *prime* if, for any ideals $J_1, J_2$, if $J = J_1 \cap J_2$, then $J = J_1$ or $J = J_2$.

**Proposition 3.2.1** *Let $N, J$ be as above.*

1. *The following are equivalent:*

    a. *$J$ is prime;*

    b. *for every $a, b \in N$, $\neg a \otimes b \in J$ or $a \otimes \neg b \in J$;*

    c. *$a \wedge b \in J$ implies $a \in J$ or $b \in J$;*

    d. *$N/J$ is totally ordered;*

e. the set of ideals extending $J$ is totally ordered by inclusion.

2. The following are equivalent:

   f. $J$ is maximal;

   g. if $a \notin J$, then there exists $k \in \mathbf{N}$ such that $(\neg a)^k \in J$;

   h. $N/J$ is isomorphic to a subalgebra of $I$. ∎

In [16, Lemma 2], Chang proves that an ideal which is maximal with respect to the property of being proper and excluding a given element of $N$, is prime. By Zorn Lemma, this implies that the intersection of all prime ideals is 0, and hence that $N$ may be embedded in $\prod\{N/p : p$ is a prime ideal of $N\}$. This is the representation of an MV-algebra as a subdirect product of totally ordered MV-algebras we were referring to in the proof of the completeness theorem. If the intersection of all maximal ideals of $N$ is 0, then by Proposition 3.2.1(2) we obtain a representation of $N$ as a subdirect product of subalgebras of $I$; we will state this in Proposition 3.2.3(3).

In [21], Di Nola shows that every MV-algebra is an algebra of "black-and-white pictures", possibly with infinitesimal grey levels:

**Theorem 3.2.2** *For every MV-algebra $N$ there exists a set $X$ and a nonstandard real unit interval $I^*$ such that $N$ is an algebra of functions $: X \to I^*$.* ∎

For any ideal $J$, let $\operatorname{Rad} J$ be the intersection of all maximal ideals extending $J$. It is known that $\operatorname{Rad} J = \{a : a \otimes (ka) \in J \text{ for every } k \geq 1\}$. We say that $N$ is:

- *simple*, if 0 is its only proper ideal;

- *semisimple*, if $\operatorname{Rad} 0 = 0$;

- *hyperarchimedean*, if all quotients of $N$ are semisimple (i.e., $\operatorname{Rad} J = J$ for every ideal $J$);

- *complete*, if $(N, \vee, \wedge, 0, 1)$ is complete as a lattice.

The element $a \in N$ is complemented (i.e., $\neg a \vee a = 1$) iff $a \oplus a = a$. As in the case of Post algebras, we denote by $\mathcal{C}(N)$ the Boolean algebra of all complemented elements of $N$.

**Theorem 3.2.3** *Let $N$ be an MV-algebra.*

1. *The following are equivalent:*

   a. *$N$ is simple;*

   b. *every nonzero element of $N$ has finite order;*

   c. *$N$ is a subalgebra of $I$.*

2. *The following are equivalent:*

   d. *$N$ is hyperarchimedean;*

   e. *prime and maximal ideals of $N$ coincide;*

   f. *for every $a \in N$, there exists $n$ such that $na \in \mathcal{C}(N)$.*

3. *The following are equivalent:*

   g. *$N$ is semisimple;*

   h. *$N$ is a subdirect product of subalgebras of $I$;*

   j. *$N$ is a subalgebra of a complete MV-algebra.*

4. *(a) implies (d), and (d) implies (g), but not conversely.*

5. *If $N$ is totally ordered, then all the conditions (a)–(j) are equivalent.*

PROOF. (1) follows from Proposition 3.2.1(2). (e) $\Longrightarrow$ (d) is clear, while (d) $\Longrightarrow$ (e) follows from the equivalence of (a) and (e) in Proposition 3.2.1(1). The equivalence of (e) and (f) is a bit scattered in the literature; a simple proof can be obtained by applying the preservation properties of the $\Gamma$ functor to [7, Théorème 14.1.2]; see also [103]. We already discussed the equivalence $(g) \Longleftrightarrow (h)$; the equivalence with (j) is proved in [50]. It is clear that (a) implies (d), which in turn implies (g). The direct product of two copies of $I$ is a nonsimple hyperarchimedean algebra, while $M_1$ is semisimple by (h), and not hyperarchimedean (Chang's algebra of Example 3.1.7 is a nonsemisimple quotient of $M_1$). Finally, assume that $N$ is totally ordered and semisimple. Then the ideals above 0 form a chain by Proposition 3.2.1(1). As $N$ is semisimple, 0 is a maximal ideal. ∎

## 3.3 Free MV-algebras

For any MV-algebra $N$, its *spectrum* $\operatorname{Spec} N$ is the set of all its prime ideals endowed with the hull-kernel topology: the closure of $U \subseteq \operatorname{Spec} N$ is $\{p \in \operatorname{Spec} N : p \supseteq \bigcap U\}$. $\operatorname{MaxSpec} N$ is the subspace of $\operatorname{Spec} N$ whose points are the maximal ideals. In what follows, $\kappa$ varies in the set $\{1, 2, 3, \ldots, \omega\}$. By [65, Proposition 8.1] the space $[0, 1]^\kappa$ is homeomorphic to $\operatorname{MaxSpec} M_\kappa$ under the mapping $\sigma \mapsto \{A \in M_\kappa : A(\sigma) = 0\}$; we freely write $A \in \sigma$ for $A(\sigma) = 0$. For $U \subseteq [0, 1]^\kappa$, we denote by $A|U$ the restriction of the function $A \in M_\kappa$ to $U$, and by $M_\kappa|U$ the MV-algebra $\{A|U : A \in M_\kappa\}$. As it was to be expected, $M_\kappa|U$ is isomorphic to $M_\kappa/\bigcap U$; all semisimple quotients of $M_\kappa$ are of this form.

We may now characterize the subsets $\Gamma$ of $FORM(L)$ for which the completeness theorem holds in the strong form:
$$C_{\operatorname{LPC}}(\Gamma) = C_\mathcal{L}(\Gamma).$$

Let $J_\Gamma$ be the ideal of $M_\omega$ generated by $\{\neg B : B \in \Gamma\}$; then $A \in C_{\operatorname{LPC}}(\Gamma)$ means $\neg A \in J_\Gamma$. On the other hand, we have

$$\begin{array}{lll}
A \in C_\mathcal{L}(\Gamma) & \text{iff} & \forall \sigma \in [0,1]^\omega (\sigma(\Gamma) \subseteq \{1\} \Rightarrow \sigma(A) = 1) \\
& \text{iff} & \forall \sigma \in [0,1]^\omega (\sigma(J_\Gamma) \subseteq \{0\} \Rightarrow \sigma(\neg A) = 0) \\
& \text{iff} & \forall \sigma \in \operatorname{MaxSpec} M_\omega (\sigma \supseteq J_\Gamma \Rightarrow \neg A \in \sigma) \\
& \text{iff} & \neg A \in \operatorname{Rad}(J_\Gamma).
\end{array}$$

So the strong completeness theorem holds exactly for those $\Gamma$'s for which $J_\Gamma = \operatorname{Rad}(J_\Gamma)$, i.e., for which $M_\omega/J_\Gamma$ is semisimple. This holds, in particular, if $\Gamma$ is finite [47, p. 84].

**Example 3.3.1** As we remarked in §2.4, adding $a \oplus a = a$ to the axioms of MV-algebras yields a subvariety of **MV** which coincides with the variety of Boolean algebras. All Boolean algebras are semisimple and hyperarchimedean, and the strong completeness theorem holds.

The $\Gamma$ functor preserves the lattice of ideals of the $\ell$-group with strong unit $(G, u)$. By well-known facts in $\ell$-group theory [7], every MV-algebra has a compact, $T_0$, and sober spectrum (*sober* means that every closed set which is not the proper union of two closed sets is the closure of a point). A difficult and still open problem is to give a topological characterization of the spectra of MV-algebras. The function which associates to every prime ideal $p$ of $N$ the unique maximal ideal that extends $p$ is a continuous retraction, and $\operatorname{MaxSpec} N$ is a compact Hausdorff space. If $N$ is hyperarchimedean, then $\operatorname{Spec} N = \operatorname{MaxSpec} N$ is totally disconnected, and coincides with the Stone space of $\mathcal{C}(N)$. In this case, $N$ is a Boolean product of subalgebras of $I$ [103].

Let $\alpha : N \to Q$ be a homomorphism of MV-algebras. $\operatorname{Spec} \alpha : \operatorname{Spec} Q \to \operatorname{Spec} N$ defined by $(\operatorname{Spec} \alpha)(p) = \alpha^{-1}(p)$ is a continuous function. If $p$ is a maximal ideal of $Q$, then $(\operatorname{MaxSpec} \alpha)(p) = \alpha^{-1}(p)$ is a maximal ideal of $N$. It follows that both Spec and MaxSpec are contravariant functors from the category of MV-algebras to the category of topological spaces and continuous functions.

**Definition 3.3.2** A *McNaughton morphism* over $[0,1]^\kappa$ is a function $f : [0,1]^\kappa \to [0,1]^\kappa$ such that, for every projection $\pi_i : [0,1]^\kappa \to [0,1]$, $f_i = \pi_i \circ f$ is a McNaughton function. If $f \circ g = g \circ f = $ (identity function), for some McNaughton morphism $g$, then we say that $f$ is a *McNaughton homeomorphism*.

If $f$ is a McNaughton morphism over $[0,1]^\kappa$, then the mapping $f^* : M_\kappa \to M_\kappa$ defined by $f^*(A) = A \circ f$ is an endomorphism of $M_\kappa$. Equivalently, $f^*$ is the endomorphism determined by $x_i \mapsto f_i$. If $\alpha : M_\kappa \to M_\kappa$ is an endomorphism, then $\operatorname{MaxSpec} \alpha$ may be naturally identified with the McNaughton morphism $f : [0,1]^\kappa \to [0,1]^\kappa$ defined by $f_i = \alpha(x_i)$. As a matter of fact, upon identifying $\sigma \in [0,1]^\kappa$ with the maximal ideal $\{A \in M_\kappa : A(\sigma) = 0\}$, we have

$$\begin{aligned}
(\text{MaxSpec}\,\alpha)(\sigma) &= \alpha^{-1}(\sigma) \\
&= \{A : \alpha(A) \in \sigma\} \\
&= \{A : (\alpha(A))(\sigma) = 0\} \\
&= \{A : (A(\ldots, f_i, \ldots))(\sigma) = 0\} \\
&= \{A : A(f(\sigma)) = 0\} \\
&= f(\sigma).
\end{aligned}$$

If we regard the monoids of endomorphisms of $M_\kappa$ and of all McNaughton functions over $[0,1]^\kappa$ as single-object categories, we trivially verify that $*$ is a contravariant isomorphism. By restriction to the invertible morphisms, we obtain an antiisomorphism between the group of automorphisms of $M_\kappa$ and the group of McNaughton homeomorphisms of $[0,1]^\kappa$.

We may now characterize McNaughton homeomorphisms. We recall that an $n$-cell is a compact convex polyhedron of affine dimension $k$. A *cellular complex* $W$ is a finite set of cells such that:

1. if $W$ contains the cell $B$, then it contains all faces of $B$;
2. every two cells of $W$ intersect in a common face.

**Theorem 3.3.3** [74] *Let $n \geq 1$, and let $f : [0,1]^n \to [0,1]^n$ be injective and continuous. Then $f$ is a McNaughton homeomorphism if and only if the following condition holds:*

*there exists a cellular complex $W$, whose $n$-cells are $B_1, \ldots, B_k$, and there exist $(n+1) \times (n+1)$ matrices $P_1, \ldots, P_k$ with integral entries, such that:*

1. $B_1 \cup \cdots \cup B_k = [0,1]^n$;
2. $P_1, \ldots, P_k$ have all determinant $+1$, or have all determinant $-1$;
3. every $P_j$ has last column of the form
$$\begin{pmatrix} 0 \\ \vdots \\ 0 \\ 1 \end{pmatrix};$$
4. $P_j$ expresses $f|B_j$ in homogeneous coordinates (i.e., if $\sigma = (s_1, \ldots, s_n) \in B_j$ and $f(\sigma) = (r_1, \ldots, r_n)$, then $(r_1, \ldots, r_n, 1) = (s_1, \ldots, s_n, 1)P_j$). ■

Let $f$ be a McNaughton homeomorphism, and let $\mu$ denote Lebesgue measure. Then, for any measurable subset $U$ of $[0,1]^n$, we have $\mu(f(U)) = \mu(U)$; this follows from the unimodularity of the matrices $P_j$.

**Example 3.3.4** The only automorphisms of $M_1$ are the identity and $x \mapsto \neg x$. In fact, the length of any segment is preserved by every McNaughton homeomorphism $f$. If $f(0) = 0$, then $f^*$ is the identity function. If $f(0) = 1$, then $f(\sigma) = 1 - \sigma$ for any $\sigma \in [0,1]$, and hence $f^* : x \mapsto \neg x$.

In addition to measure, McNaughton homeomorphisms preserve denominators. Let $\sigma = (s_1, \ldots, s_n)$ be a point of $[0,1]^n$ having rational coordinates. There exist uniquely determined relatively prime integers $b_1, \ldots, b_{n+1}$ such that $b_{n+1} \geq 1$ and $s_i = b_i/b_{n+1}$, for every $i = 1, \ldots, n$. We say that $b_{n+1}$ is the *denominator* of $\sigma$, and we define the $k$-grid $Gr_k$ of $[0,1]^n$ to be the set of all rational points in $[0,1]^n$ having denominator $k$. Every McNaughton homeomorphism induces, for any $k$, a permutation of the $k$-grid onto itself. Let $G$ be the group of McNaughton homeomorphisms over $[0,1]^n$. For every $k$, we have a natural homomorphism : $G \to \text{Sym}(Gr_k)$; let $H_k$ be the kernel. By continuity, the intersection of all the $H_k$'s is the identity function, so that $G$ is residually finite. Let $I_{k+1}$ be the subalgebra of $I$ generated by $1/k$. Then $M_n|Gr_k$ is the free MV-algebra over $n$ generators in the subvariety of **MV** generated by $I_{k+1}$, and preservation of denominators means simply that every automorphism of $M_n$ induces an automorphism of $M_n|Gr_k$. It is not the case that every automorphism of $M_n|Gr_k$ extends to an automorphism of $M_n$; this follows easily from Example 3.3.4, by taking $n = 1$ and $k = 3$.

# 4 Why many-valued logic?

We now try to sum up our discussions, and to recapitulate the distinguishing features of many-valued logic.

First of all, we think that many-valued logic is a mathematically interesting theory, with relatively deep theorems, and connections with other branches of mathematics. Some of these connections are apparent in §3; links with lattice theory and universal algebra are everywhere, while the classification of connectives involves functional equations.

Secondly, certain issues of classical logic are clarified by the many-valued extensions. Such issues include the rôle of quantifiers, the introduction-elimination rules for connectives, the relations among the various proof systems, the properties of consequence operations, the rôle of negation.

Concerning negation, we are now forced to distinguish various cases of inconsistency/contradiction for a set of formulas $\Gamma$:

1. the case in which every formula is deducible from $\Gamma$;

2. the case in which $\Gamma$ proves a formula $A$ along with its negation $\neg A$;

3. the case in which $\Gamma$ proves a specific formula (or set of formulas) which is assumed to express falsehood.

In particular, we can have a logic in which contradiction does not imply inconsistency, i.e., in which (2) holds, but not (1). Such logics are, e.g., Belnap's logic, or Mundici's formalization of the Ulam game with lies. The introduction of new truth-values has also been proposed as an alternative to negation-as-failure in logic programming [20]. We consider this to be a third remarkable feature.

As a fourth point, we stress that many-valued logic allows one to work truth-functionally on imprecise statements, provided complete information is given. Recall the distinction between degrees of truth ("The sky is rather blue") and degrees of belief ("We are confident that tomorrow will not rain"). The former statement is not affected by increase of information, while the latter is. Although belief and probability are not truth-functional, truth can be treated truth-functionally. Once this distinction is clear, we may be satisfied with Scott's interpretation of the Łukasiewicz connectives, provided we read *degree of falsehood* for *degree of error*. Have a look again at the implication: $A \to B$ means now naturally "$B$ is at least as true as $A$". Hence the degree of falsehood of $A \to B$ in a state of affairs $\sigma$ is just the amount by which $\sigma(A)$ dominates $\sigma(B)$: in symbols

$$\sigma(A \to B) = \begin{cases} 1, & \text{if } \sigma(A) \leq \sigma(B); \\ 1 - (\sigma(A) - \sigma(B)), & \text{otherwise}; \end{cases}$$

which is Łukasiewicz implication. We already saw that the resulting non-idempotence of $\oplus$ and $\otimes$ may be a welcome feature. Note that we are not claiming that we know what a degree of truth is, nor that we must treat them truth-functionally. We are merely saying that there are statements to which it is reasonable to attach a label different from $\{true, false\}$, and that sometimes we can work out propagation rules for these labels.

Let us try an example: consider an art critic, expert in Pietro della Francesca's works. He has detailed information about Pietro's paintings, and is fully confident with the artist style. A new painting is discovered, whose attribution is uncertain. How does the critic work? We can formalize the situation by declaring that the critic's knowledge is contained in a database $\Omega$, that examining the new picture he forms a second database $\omega$, and that he decides the attribution by confronting $\Omega$ and $\omega$. The kind of information contained in the two databases is hardly numeric or boolean. Rather, it is expressed in the form of statements like "In these circumstances, Pietro tends to use that colour" or "The posture of the hands is typical". This kind of statements are best visualized as fuzzy predicates over a set, as defined in §2.6. The comparison results in an averaging process, and its output is a truth value. Note that a numeric answer between 0 and 1 is reasonable, and may even be unavoidable: indeed, many paintings in Pietro's times are works of the master's school, the artist painting the main figures only. Note also that some of the information contained in $\Omega$ may result from works of dubious attribution. If so, $\Omega$ is affected by a further level of approximation, and the final answer will be something like "It is plausible that this work is mostly due to Pietro". In next Chapter's notation, we have the evaluated formula $[mostly\_authentic(picture); plausible]$.

This extension of our example leads to our final point: many-valued logic provides the theoretical basis to fuzzy logic in the narrow sense. The picture we have in mind is the following:

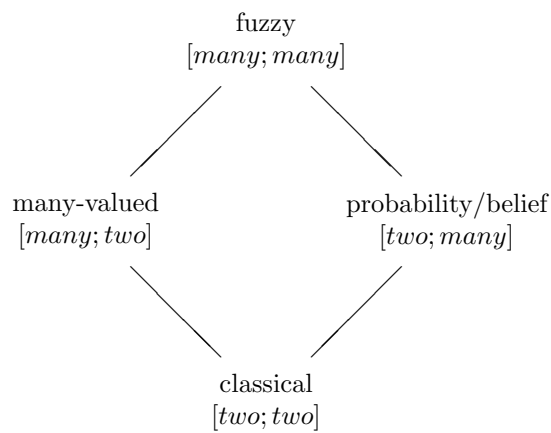

It is suggestive to think of fuzzy logic as a many-valued logic in which the imprecision migrates from facts to statements about facts, i.e., from the language to the metalanguage; in the evaluated formula notation, from the left to the right-hand side of ";". People working in fuzzy logic may benefit from the by now well established corpus of techniques of many-valued logics; to cite an example, the original Pavelka proof [77] of the completeness theorem for fuzzy logic has been drastically simplified in [43] by making use of the classical completeness theorem of Łukasiewicz logic.